\documentclass[a4paper,12pt]{article}

\usepackage{amsmath, amssymb, amsthm}

\usepackage{color}
\definecolor{vio}{rgb}{0.5,0,0.5}
\definecolor{gre}{rgb}{0.1,0.6,0}

\usepackage{titlesec}
\titleformat{\section}{\bfseries}{\thesection}{1em}{}
\titleformat{\subsection}{\itshape}{\thesubsection}{1em}{}

\numberwithin{equation}{section}

\def\expe{\hbox{\rm e}}
\def\ve{\varepsilon}
\def\vp{\varphi}
\def\hate{\hat\theta}
\def\dd{\,\mathrm{d}}
\def\dive{\mathrm{\,div\,}}

\def\dch{(\chi{+}\rho^*(1{-}\chi))}

\def\eau{H\!${}_2$\!O }
\def\on{^{(n)}}
\def\om{^{(m)}}
\def\oi{^{(i)}}

\def\real{\mathbb{R}}
\def\nat{\mathbb{N}}

\def\supess{\mathop{{\rm sup\,ess}\,}}

\def\nas{\nabla_s}
\def\iv{\int_V}
\def\ipv{\int_{\partial V}}
\def\io{\int_{\Omega}}
\def\ipo{\int_{\partial\Omega}}

\def\for{\hbox{\ for}\ }

\newfont{\ctv}{msam10}

\newcommand{\bbox}{\mbox{\ctv \symbol{4}}}
\def\QED{{${}\hfill\bbox$}}
\newenvironment{pf}[1]{\par\vskip1mm{\noindent\it #1.}\ }{\QED\par
\vskip2mm}

\def\bpf{\begin{pf}}
\def\epf{\end{pf}}

\def\be{\begin{equation}\label}
\def\ee{\end{equation}}

\def\ber{\begin{eqnarray}}
\def\eer{\end{eqnarray}}

\def\bers{\begin{eqnarray*}}
\def\eers{\end{eqnarray*}}

\def\bpf{\begin{pf}}
\def\epf{\end{pf}}

\newtheorem{theorem}{Theorem}[section]
\newtheorem{lemma}[theorem]{Lemma}

\newtheorem{hypothesis}[theorem]{Hypothesis}
\newtheorem{proposition}[theorem]{Proposition}

\begin{document}

\title{Unsaturated deformable porous media flow \\
with phase transition\footnote{The financial supports of the FP7-IDEAS-ERC-StG \#256872
(EntroPhase), of the project Fondazione Cariplo-Regione Lombardia  MEGAsTAR 
``Matematica d'Eccellenza in biologia ed ingegneria come accelleratore di una nuona strateGia per l'ATtRattivit\`a dell'ateneo pavese'', and of the GA\v CR Grant GA15-12227S and RVO: 67985840
are gratefully acknowledged. The paper 
also benefited from the support the GNAMPA (Gruppo Nazionale per l'Analisi Matematica, la Probabilit\`a e le loro Applicazioni) of INdAM (Istituto Nazionale di Alta Matematica) for ER.}}

\author{Pavel Krej\v c\'{\i}
\thanks{Institute of Mathematics, Czech Academy of Sciences, \v Zitn\'a~25,
CZ-11567~Praha 1, Czech Republic, E-mail {\tt krejci@math.cas.cz}.}\,\,,
Elisabetta Rocca
\thanks{Dipartimento di Matematica, Universit\`a degli Studi di Pavia and IMATI-C.N.R., Via Ferrata 5, 
I-27100 Pavia, Italy, E-mail {\tt elisabetta.rocca@unipv.it}.}
\,, and J\"urgen Sprekels
\thanks{Weierstrass Institute for Applied
Analysis and Stochastics, Mohrenstrasse~39, D-10117 Berlin,
Germany, E-mail {\tt  sprekels@wias-berlin.de}, and Department of Mathematics, Humboldt-Universit\"at
zu Berlin, Unter den Linden 6, D-10099 Berlin, Germany.}
}

\maketitle

\begin{abstract}\noindent
In the present paper, a continuum model is introduced for fluid flow in a 
deformable porous medium, where the fluid may undergo phase transitions. Typically, such problems arise
in modeling liquid-solid phase transformations in groundwater flows.
The system of equations is derived here from the conservation principles for mass, momentum,
and energy and from the Clausius-Duhem inequality for entropy. It
couples the evolution of the displacement in the matrix material, of the capillary pressure,
of the absolute temperature, and of the phase fraction.
Mathematical results are proved under the additional hypothesis that
inertia effects and shear stresses can be neglected. For the resulting
highly nonlinear system of two PDEs, one ODE and one ordinary differential inclusion with
natural initial and boundary conditions, existence of global in time solutions is proved
by means of cut-off techniques and suitable Moser-type estimates.
\end{abstract}


\section*{Introduction}\label{int}

A model for fluid flow in partially saturated porous media with thermomechanical interaction
was proposed and analyzed in \cite{ak,dkr1, dkr2}.
Here, we extend the model by including the effects of freezing and melting of the fluid in the pores.
Typical examples, in which such situations arise, are related to groundwater flows and to
the freezing-melting cycles of water sucked into the pores of concrete. Notice that the latter
process forms one of the main 
reasons for the degradation of concrete in buildings, bridges, and roads. However, many of the governing
effects in concrete like the multi-component microstructure, the breaking of pores, chemical reactions,
the hysteresis of the saturation-pressure curves, and the occurrence of shear stresses, are still
neglected in our model.    

While often in continuum models for three-component and multi-component porous media the
intention is to describe the
propagation of sound waves in these media (e.g., \cite{albers,wil}), we investigate in the present paper -- instead of using partial balance equations for each component -- a continuum model 
combining the principles of conservation of mass and momentum with
the first and the second principles of thermodynamics. In, e.~g., \cite{albers}, flow in
porous media is described in Eulerian coordinates in order to incorporate, for example,
the effects of fast convection. Here, instead, we assume that slow diffusion is dominant,
and choose the Lagrangian description as in \cite{ak, dkr1, wil}. The resulting system of
coupled ODEs and PDEs  then appears to be a nonlinear
extension of the linear model in \cite{show}, referred to as a simplified Biot system,
to the case when also the occurrence of temperature changes and
phase transitions is taken into account. In addition to the model studied in \cite{dkr1},
we include here the effects of freezing and melting. 
The idea is the following. The pores in the matrix material contain a mixture of \eau and gas,
and \eau itself is a mixture of the liquid (water) and the solid phase (ice).
That is, in addition to the other physical quantities like capillary pressure, displacement, and
absolute temperature, we need to consider the evolution of a phase parameter
$\chi$ representing the relative proportion of water in the \eau part and its influence on pressure changes due to the different
mass densities of water and ice. Unlike in \cite{ak,dkr1,dkr2,ss}, we do not consider hysteresis
in the model. We believe that the mathematical results can be extended to the case of capillary
hysteresis as in \cite{ak,dkr1,dkr2}. In our model without shear stresses, elastoplastic hysteresis
effects as in \cite{ak,dkr1,ss} cannot occur.

As it will be detailed in Section~\ref{mod}, we assume that the deformations are small, so that $\dive u$ is the relative local volume change, where $u$ represents the displacement vector. Moreover, we assume that the volume of the matrix material does not change during the process, and thus the volume and mass balance equations
with Darcy's law for the water flux lead to a nonlinear degenerate parabolic equation for the
capillary pressure, see~\eqref{e4b}. In the equation of motion, we take into account the pressure components due to phase transition and temperature changes, and we further simplify the system in order to make it mathematically tractable by assuming that the process is quasistatic and the shear stresses are negligible. The problem of existence of solutions for the coupled system without this assumption is open and, in our opinion, very challenging. Finally, we use the balance of internal energy and the entropy inequality to derive the 
dynamics for absolute temperature and phases; they turn out to be, respectively, a parabolic equation for the 
temperature with highly nonlinear right-hand side (quadratic in the derivatives) and an ordinary 
differential inclusion for the phase parameter $\chi$. 

Finally, let us note that -- in order to model the freezing and melting phenomena 
in the pores -- we have borrowed here some ideas from our earlier publications
on freezing and melting in containers filled with water with rigid, elastic, or elastoplastic boundaries
(cf. \cite{mb,kr,krsbottle,krsgrav,krsWil}). It was shown there how important it is to account for
the difference in specific volumes of water and of ice.

There is an abundant classical mathematical literature on phase
transition processes, see, e.g., the monographs \cite{bs},
\cite{fremond}, \cite{visintin}, and the references therein. It
seems, however, that only few publications take into account
that the mass densities and specific volumes of the phases differ. In
\cite{fr1}, the authors proposed to interpret a phase transition
process in terms of a balance equation for macroscopic motions,
and to include the possibility of voids. Well-posedness of an
initial-boundary value problem associated with the resulting PDE
system was proved there, and the case of two different densities
$\varrho_1$ and $\varrho_2$ for the two substances undergoing
phase transitions was pursued in \cite{fr2}.

Let us also mention the papers \cite{rr1, rr2, RocRos12, RocRos14} dealing
with macroscopic stresses in phase transitions models, where the
different properties of the viscous (liquid) and elastic (solid)
phases were taken into account and the coexisting viscous and
elastic properties of the system were given a distinguished role,
under the working assumption that they indeed influence the phase
transition process. The model studied there includes inertia,
viscous, and shear viscosity effects (depending on the phases). This is reflected in the analytical
expressions of the associated PDEs for the strain $u$ and
the phase parameter $\chi$: the $\chi$-dependence, e.~g.,  in the
stress-strain relation, leads to the possible degeneracy of the
elliptic operator therein. Finally, we can quote in this framework the
model analyzed  in~\cite{kss2} and \cite{kss}, which
pertains to nonlinear thermoviscoplasticity: in the
spatially one-dimensional case, the authors prove the global
well-posedness of a PDE system that both incorporates hysteresis
effects and models phase change but, however, does not display
a degenerating character.

Another coupled system for temperature, displacement, and phase parameter
has been derived in order to model the full thermomechanical behavior of
shape memory alloys. A long list of references for further developments
can be found in the monographs \cite{fremond} and \cite{visintin}.

The paper is organized as follows: in the next Section~\ref{mod}, we derive the model
in full generality from the basic principles of continuum thermodynamics.
In Section~\ref{mat}, we state the mathematical problem, the main assumptions on the data,
and the main Theorem \ref{t1}, 
the proof of which is split into Sections~\ref{cut}, \ref{apr}, and \ref{proo}. The
steps of the proof are as follows: we first cut off some of the pressure and temperature dependent terms in the system in Section~\ref{cut} by means of a cut-off parameter $R$ and solve the related problem employing a special Galerkin approximation scheme. Then, in Section~\ref{apr}, we first prove the positivity of the temperature by means of a maximum principle technique, and then we perform the -- independent of  $R$ --  estimates on the system. They mainly consist of: the energy estimate, the so-called Dafermos estimate (with negative small powers of the temperature), Moser-type and then higher-order estimates for the capillary pressure and for the temperature. This allows us in Section~\ref{proo} to pass to the limit
in the cut-off system as $R\to\infty$, which will conclude the proof of the existence result. 


\section{The model}\label{mod}
We consider a connected domain $\Omega\subset \real^3$ filled by a deformable matrix material with pores
containing a mixture of \eau and gas, where we assume that \eau may appear in one of the two phases:
water or  ice. We also assume that the volume of the solid matrix remains constant during the process,
and let $c_s \in (0,1)$  be the relative proportion of solid in the total reference volume.
We denote, for $x\in \Omega$ and time $t\in [0,T]$,
\begin{description}
\item $W(x,t)\in [0,1]$ ... relative proportion of \eau in the total pore volume;
\item $A(x,t)\in [0,1]$ ... relative proportion of gas in the total pore volume;
\item $\chi(x,t)\in [0,1]$ ... relative proportion of water in the \eau part;
\item $\xi(x,t)$ ... mass flux vector;
\item $p(x,t)$ ... capillary pressure;
\item $u(x,t)$ ... displacement vector;
\item $\sigma(x,t)$ ... stress tensor;
\item $\theta(x,t)$ ... absolute temperature.
\end{description}
Then $\chi W$ represents the relative proportion of water in the total pore volume, and
$(1-\chi) W$ represents the relative proportion of ice in the total pore volume.

We assume that the deformations are small, so that $\dive u$ is the relative local volume change.
By hypothesis, the volume of the matrix material does not change, so that the volume balance reads
\be{e1}
W(x,t) + A(x,t) + c_s = 1+\dive u(x,t).
\ee

For $A$, we assume the functional relation
\be{e2}
A = 1-c_s - \vp(p)\,,
\ee
where $\vp$ is an increasing function that satisfies $\vp(-\infty) = \vp^\flat \in (0,1)$
and $\vp(\infty) = 1-c_s$, $\vp^\flat + c_s < 1$.
This means that the porous medium cannot be made completely dry by thermomechanical processes alone.
Combining \eqref{e1} with \eqref{e2}, we obtain that
\be{e1a}
W = \vp(p) + \dive u\,.
\ee

\subsection{Mass balance}\label{mass}

Consider an arbitrary control volume $V \subset \Omega$. The water content in $V$ is given
by the integral $\iv \rho_L \chi W \dd x$, where $\rho_L$ is the water mass density, and
the ice content is $\iv \rho_S (1-\chi) W \dd x$, where $\rho_S$
is the ice mass density. The mass conservation principle then reads
\be{e3}
\frac{\dd}{\dd t} \iv \rho_L \chi W \dd x + \ipv \xi\cdot n\dd s(x)
= -\frac{\dd}{\dd t} \iv \rho_S  (1-\chi) W \dd x\,,
\ee
where $n$ the unit outward normal vector to $\partial V$. In differential form, we obtain
\be{e4}
\rho_L (\chi W)_t + \dive \xi = -\rho_S  ((1-\chi) W)_t\,.
\ee
The right-hand side of \eqref{e4}
is the positive or negative liquid water source due to the solidification or melting of the ice.
We assume the water flux in the form of the Darcy law
\be{e5}
\xi = - \mu(p) \nabla p,
\ee
with a proportionality factor $\mu(p)>0$. This, \eqref{e1a}, and \eqref{e4}, yield the equation
\be{e4b}
\big(\dch(\vp(p) + \dive u)\big)_t - \frac{1}{\rho_L}\dive (\mu(p)\nabla p) = 0,
\ee
with $\rho^* = \rho_S/\rho_L \in (0,1)$.


\subsection{Equation of motion}\label{moti}

The equation of motion is considered in the form
\be{mo1}
\rho_M u_{tt} - \dive \sigma = g\,,
\ee
where $\rho_M$ is the mass density of the matrix material, $\sigma$ is the stress tensor,
and $g$ is a volume force acting on the body (e.~g., gravity). For $\sigma$, we prescribe
the constitutive equation
\be{mo2}
\sigma = B\ve_t + A\ve + \big(\dch(\lambda \dive u - p) - \beta(\theta-\theta_c)\big)\delta\,,
\ee
where $\ve = \nas u :=\frac12 (\nabla u + \nabla u^T)$ is the small strain tensor, $\delta$
is the Kronecker tensor, 
$B$ is a symmetric positive definite viscosity tensor, $A$
is the symmetric positive definite elasticity tensor of the matrix material, $\lambda>0$ is the
bulk elasticity modulus of water,
$\theta>0$ is the absolute temperature, $\theta_c>0$ is a fixed referential temperature,
and $\beta \in \real$ is the relative solid-liquid thermal expansion coefficient. The term
$\,\dch(\lambda \dive u - p)\,$ accounts for the pressure component
due to the phase transition.


\subsection{Energy and entropy balance}\label{ener}

We have to derive formulas for the densities of internal energy $U$ and entropy $S$
such that the energy balance balance equation and the Clausius--Duhem inequality hold for all processes.
Let $q$ be the heat flux vector, and let $V \subset \Omega$ be again an arbitrary control volume.
The total internal energy in $V$ is $\iv U \dd x$, and the total mechanical power $Q(V)$ supplied to $V$
equals
$$
Q(V) = \iv \sigma:\ve_t \dd x - \ipv \frac{1}{\rho_L} p\, \xi\cdot n \dd s(x)\,,
$$
where $\xi$ is the fluid mass flux \eqref{e5}. We thus have that
\be{e11}
\frac{\dd}{\dd t} \iv U\dd x + \ipv q\cdot n \dd s(x) =
 \iv \sigma:\ve_t \dd x - \ipv \frac{1}{\rho_L} p\,\xi \cdot n \dd s(x)\,.
\ee
Again, by the Gauss formula, we obtain the energy balance equation in differential form, namely
\be{e12}
U_t + \dive q = \sigma:\ve_t - \frac{1}{\rho_L} \dive( p\xi)\,.
\ee
The internal energy and entropy densities $U$ and $S$, as well as the heat flux vector $q$,
have to be chosen in order to satisfy, for all processes, the Clausius--Duhem inequality
\be{e13}
S_t + \dive\big(\frac{q}{\theta}\big) \ge 0,
\ee
or, taking into account the energy balance \eqref{e12},
\be{e14}
U_t - \theta S_t + \frac{q\cdot\nabla \theta}{\theta} \le \sigma:\ve_t - \frac{1}{\rho_L} \dive( p\xi)\,.
\ee
We consider $\ve, \chi, p, \theta$ as state variables and $U, S$ as state functions, independent of $\nabla\theta$.
Hence, as a consequence of \eqref{e14}, two inequalities have to hold separately for all processes, namely
\be{e15}
q\cdot\nabla \theta \le 0\,, \quad U_t - \theta S_t \le \sigma:\ve_t - \frac{1}{\rho_L} \dive( p\xi)\,.
\ee
For simplicity, we assume Fourier's law for the heat flux,
\be{e16}
q = -\kappa(\theta) \nabla\theta,
\ee
with the heat conductivity coefficient $\kappa = \kappa(\theta) >0$. 
We further introduce the free energy $F$ by the formula
$F = U - \theta S$, so that, in terms of $F$, the second inequality in \eqref{e15} takes the form
\be{e17}
\quad F_t + \theta_t S \le \sigma:\ve_t + \frac{1}{\rho_L} \dive( p\mu(p)\nabla p)\,.
\ee
We claim that the right choice of $F$ for \eqref{e17} to hold is given by
\ber \nonumber
F &=& \frac12 A\ve:\ve + \dch\left(V(p)+ \frac{\lambda}{2} (\dive u)^2\right)\\ \label{e18f}
 && +\, L\chi \left(1 - \frac{\theta}{\theta_c}\right)
 -\beta(\theta - \theta_c)\dive u  + F_0(\theta) + I(\chi),\\ \label{e18s}
S &=& -\frac{\partial F}{\partial\theta} = \frac{L}{\theta_c}\chi + \beta \dive u - F'_0(\theta),
\eer
where
\be{e18a}
V(p) = p \vp(p) - \Phi(p)\,, \quad \Phi(p) = \int_0^p \vp(\tau) \dd \tau\,,
\ee
$F_0(\theta)$ is a purely caloric component of $F$, $L>0$ is the latent heat, and $I$
is the indicator function of the interval $[0,1]$. It remains to check that if we
choose the phase dynamics equation in the form
\be{e6}
\gamma(\theta) \chi_t+ \partial I(\chi)
\ni (1-\rho^*)\left(\Phi(p)+ p\dive u  - \frac{\lambda}{2}(\dive u)^2\right)
+ L\left(\frac{\theta}{\theta_c}-1\right)
\ee
with a coefficient $\gamma(\theta)>0$, then \eqref{e17} holds for all processes.
Indeed, by \eqref{e18f}--\eqref{e18s} and \eqref{e4b} we have that
\bers
F_t + \theta_t S &=& A\ve: \ve_t + \dch(V'(p)p_t + \lambda \dive u \dive u_t)\\
&&+\, (1 - \rho^*)\chi_t \left(V(p)+ \frac{\lambda}{2} (\dive u)^2\right)
+ L\chi_t \left(1 - \frac{\theta}{\theta_c}\right)-  \beta(\theta-\theta_c)\dive u_t\,,\\
\sigma:\ve_t &=& B \ve_t:\ve_t + A\ve: \ve_t + \dch(\lambda \dive u \dive u_t - p \dive u_t)\\
&& -\, \beta(\theta-\theta_c)\dive u_t\,,\\
 \frac{1}{\rho_L}\dive (p \mu(p)\nabla p) &=& \frac{1}{\rho_L} \mu(p)|\nabla p|^2 +
p \dch(\vp'(p)p_t + \dive u_t)\\
&& +\, p(1 - \rho^*)\chi_t (\vp(p) + \dive u)\,.
\eers
Hence (note that $V(p) - p\vp(p) = -\Phi(p)$),
\ber \nonumber
F_t + \theta_t S - \sigma:\ve_t - \frac{1}{\rho_L}\dive (p \mu(p)\nabla p) &=&
 - B \ve_t:\ve_t - \frac{1}{\rho_L} \mu(p)|\nabla p|^2 \\ \nonumber
 && \hspace{-25mm} +\, \chi_t \left(L\left(1 - \frac{\theta}{\theta_c}\right)
 + (1 - \rho^*)\left( \frac{\lambda}{2} (\dive u)^2 - \Phi(p) - p \dive u\right)\right)\\ \label{xx}
 &=& - B \ve_t:\ve_t - \frac{1}{\rho_L} \mu(p)|\nabla p|^2 - \gamma(\theta)\chi_t^2 \le 0,
\eer
by virtue of \eqref{e6}, so that \eqref{e17} holds. 

Now observe that
\ber \nonumber
U &=& F + \theta S\\ \nonumber
&=&
\frac12 A\ve:\ve + \dch\left(V(p)+ \frac{\lambda}{2} (\dive u)^2\right)\\ \label{e20}
&&+\, L\chi + \beta\theta_c\dive u  + F_0(\theta)-\theta F_0'(\theta) + I(\chi)\,.
\eer
The derivative of the purely caloric component $F_0(\theta) - \theta F'_0(\theta)$ is the specific heat
capacity $c(\theta) = -\theta F''(\theta)$. Assuming that $c(\theta) = c_0$ is a positive constant,
we obtain that $F_0(\theta) = - c_0\theta \log(\theta/\theta_c)$ up to a linear function, and
\be{e20a}
U = \frac12 A\ve:\ve + \dch\left(V(p)+ \frac{\lambda}{2} (\dive u)^2\right)
+ L\chi + \beta\theta_c\dive u + c_0 \theta + I(\chi)\,.
\ee
We now rewrite Eq.~\eqref{e12} in a more suitable form, using \eqref{xx}. We have
\ber \nonumber
0 &=& U_t + \dive q - \sigma:\ve_t - \frac{1}{\rho_L} \dive( p\mu(p)\nabla p)\\ \nonumber
&=& (F+\theta S)_t + \dive q - \sigma:\ve_t - \frac{1}{\rho_L} \dive( p\mu(p)\nabla p)\\ \label{ene4}
&=& - B \ve_t:\ve_t - \frac{1}{\rho_L} \mu(p)|\nabla p|^2 - \gamma(\theta)\chi_t^2 +
\theta S_t + \dive q\,,
\eer
which yields the identity
\be{ene5}
c_0 \theta_t - \dive(\kappa(\theta)\nabla\theta) =  B \ve_t:\ve_t
+\frac{1}{\rho_L} \mu(p)|\nabla p|^2 + \gamma(\theta)\chi_t^2 - \frac{L}{\theta_c}\theta\chi_t
- \beta \theta \dive u_t\,.
\ee


\section{The mathematical problem}\label{mat}

We consider the system
\begin{align}\label{ae1}
&\big(\dch(\vp(p) + \dive u)\big)_t \,=\, \frac{1}{\rho_L}\dive (\mu(p)\nabla p)\,,\\ 
\label{ae2}
&\,\rho_M u_{tt} \,=\, \dive \sigma + g\,,\\[2mm] 
&\hspace*{8.7mm}\sigma \,=\, B\nas u_t + A\nas u + (\dch(\lambda \dive u - p)
\,-\, \beta(\theta-\theta_c))\delta\,,\label{ae3} 
\end{align}
\begin{align}
\label{ae4}
&\hspace*{8.5mm}\gamma(\theta) \chi_t+ \partial I(\chi)
\,\ni\, (1-\rho^*)\left(\Phi(p)+ p\dive u  - \frac{\lambda}{2}(\dive u)^2\right)
\,+\, L\left(\frac{\theta}{\theta_c}-1\right),\\[2mm]
&\label{ae5}
c_0 \theta_t - \dive(\kappa(\theta)\nabla\theta)\,=\, B \nas u_t:\nas u_t
+\frac{1}{\rho_L} \mu(p)|\nabla p|^2 + \gamma(\theta)\chi_t^2
\,-\, \frac{L}{\theta_c}\theta\chi_t - \beta \theta \dive u_t\,,
\end{align}
for the unknown functions $p,u,\chi,\theta$, coupled with the boundary conditions
\ber\label{be1}
u &\!\!=\!\!& 0\,,\\ \label{be2}
\xi\cdot n &\!\!=\!\!& \alpha(x) (p - p^*)\,,\\ \label{be3}
q \cdot n &\!\!=\!\!& \omega(x) (\theta - \theta^*)\,,
\eer
on $\partial\Omega$, where $p^*$ is a given outer pressure, $\theta^*$ is a given outer temperature,
$\alpha(x) \ge 0$ is the permeability of the boundary, and $\omega(x) \ge 0$
is the heat conductivity of the boundary.

We can also simplify the problem by assuming that water is incompressible.
This corresponds to the choice $\lambda = 0$, whence the system becomes
\ber\label{le1}
\big(\dch(\vp(p) + \dive u)\big)_t &\!\!=\!\!& \frac{1}{\rho_L}\dive (\mu(p)\nabla p)\,,
\\ \label{le2}
\rho_M u_{tt} &\!\!=\!\!& \dive \sigma + g\,,\\  \label{le3}
\sigma &\!\!=\!\!& B\nas u_t + A\nas u - (p \dch - \beta(\theta-\theta_c))\delta\,,\qquad\\ \label{le4}
\gamma(\theta) \chi_t+ \partial I(\chi)
&\!\!\ni\!\!& (1-\rho^*)\left(\Phi(p)+ p\dive u\right) + L\left(\frac{\theta}{\theta_c}-1\right),\\ \nonumber
c_0 \theta_t - \dive(\kappa(\theta)\nabla\theta) &\!\!=\!\!&  B \nas u_t:\nas u_t
+\frac{1}{\rho_L} \mu(p)|\nabla p|^2 + \gamma(\theta)\chi_t^2\\ \label{le5}
&&-\, \frac{L}{\theta_c}\theta\chi_t - \beta \theta \dive u_t\,.
\eer
We further simplify the system by assuming that the process is quasistatic and that the shear stresses
are negligible. Then \eqref{le2}--\eqref{le3} can be reduced to 
\ber \label{le2a}
0 &\!\!=\!\!& \dive \sigma + g\,,\\  \label{le3a}
\sigma &\!\!=\!\!& (\nu \dive u_t + \lambda_M \dive u - p \dch - \beta(\theta-\theta_c))\delta\,.
\eer
Assuming that the force $g$ admits a potential $G$, that is, $g = \nabla G$, this yields
\be{le6}
\nu \dive u_t + \lambda_M \dive u - p \dch - \beta(\theta-\theta_c) \,=\, -G + H(t)\,,
\ee
where $H(t)$ is an ``integration constant'', $\nu$ is the bulk viscosity coefficient, and
$\lambda_M$ is the bulk elasticity modulus of the matrix material. In view of the boundary
condition \eqref{be1}, we have that
\be{le7}
H(t)\, =\, - \frac{1}{|\Omega|}\io (p \dch + \beta(\theta-\theta_c) -G)(x,t)\dd x\,.
\ee
With the new unknown function $w = \dive u$, which represents the {\em relative volume change\/},
the system \eqref{le1}--\eqref{le5} then becomes

\ber\label{lu1}
\big(\dch(\vp(p) + w)\big)_t &\!\!=\!\!& \frac{1}{\rho_L}\dive (\mu(p)\nabla p)\,,
\\[1mm]
\label{lu2}
\nu w_t + \lambda_M w &\!\!=\!\!&  p \dch + \beta(\theta-\theta_c) - G + H(t)\,,\\[1mm]
\label{lu4}
\gamma(\theta) \chi_t+ \partial I(\chi)
&\!\!\ni\!\!& (1-\rho^*)\left(\Phi(p)+ p w\right) + L\left(\frac{\theta}{\theta_c}-1\right),\\[1mm]
\label{lu5}
c_0 \theta_t - \dive(\kappa(\theta)\nabla\theta) &\!\!=\!\!&  \nu w_t^2
+\frac{1}{\rho_L} \mu(p)|\nabla p|^2 + \gamma(\theta)\chi_t^2
- \frac{L}{\theta_c}\theta\chi_t - \beta \theta w_t\,.
\eer
We prescribe the initial conditions
\begin{align} \label{inip}
p(x,0) &\,=\, p^0(x)\,,\\ \label{iniu}
w(x,0) &\,=\, w^0(x)\,,\\ \label{inich}
\chi(x,0) &\,=\, \chi^0(x)\,,\\ \label{init}
\theta(x,0) &\,=\, \theta^0(x)\,.
\end{align}
The weak formulation of Problem \eqref{lu1}--\eqref{lu5} reads as follows:
\ber\label{wu1}
\io\left((\dch (\vp(p)+w))_t \eta +\frac{1}{\rho_L}\mu(p)\nabla p\cdot\nabla \eta\right)\dd x
&\!\!=\!\!& \ipo \alpha(x)(p^*-p)\eta \dd s(x),\qquad \\ \label{wu2}
\nu w_t + \lambda_M w  - p \dch - \beta(\theta-\theta_c) &\!\!=\!\!& - G + H(t)\quad\mbox{a.~e.},\\ \label{wu4}
 \gamma(\theta) \chi_t + \partial I(\chi) - (1-\rho^*)(\Phi(p)+pw) &\!\!\ni\!\!& L\left(\frac{\theta}{\theta_c}-1\right)
\quad\mbox{a.~e.},\\ \label{wu5}
\io\left(c_0\theta_t - \gamma(\theta)\chi_t^2
+ \frac{L}{\theta_c}\theta\chi_t -  \nu w_t^2 + \beta \theta w_t\right) \zeta\dd x && \\ \nonumber
+\io\left(- \frac{1}{\rho_L}\mu(p)|\nabla p|^2 \zeta + \kappa(\theta) \nabla \theta\cdot\nabla \zeta\right)\dd x
&=& \ipo \omega(x)(\theta^*-\theta)\zeta \dd s(x)\,,
\eer
almost everywhere in $(0,T)$ and for all test functions $\eta\in W^{1,2}(\Omega)$ and $\zeta \in W^{1,q^*}(\Omega)$,
with some $q^* > 1$ that will be specified below in Theorem \ref{t1}.

\vspace*{3mm}
\begin{hypothesis}\label{h1} 
{\rm We fix a time interval $[0,T]$ and assume that the data of Problem \eqref{wu1}--\eqref{wu5} have the following
properties:}
\begin{itemize}
\item[{\rm (i)}] $\gamma : [0,\infty) \to [0,\infty)$ {\rm is continuous;
$\exists 0 < c_\gamma < C_\gamma: c_\gamma(1+\theta) \le \gamma(\theta) \le C_\gamma(1+\theta)$ for all $\theta \ge 0$;}
\item[{\rm (ii)}] $\kappa : [0,\infty) \to [0,\infty)$ {\rm is continuous; $\exists 0<c_\kappa< C_\kappa$,
$0< a < 1$, $a< \hat a < \frac{16}{5} + \frac65 a:
c_\kappa (1+ \theta^{1+a}) \le \kappa(\theta) \le C_\kappa (1+ \theta^{1+\hat a})$ for all $\theta\ge 0$;}
\item[{\rm (iii)}] {\rm $\theta^0 \in W^{1,2}(\Omega) \cap L^\infty(\Omega)$,
$\theta^*\in L^\infty(\partial\Omega\times (0,T))$, $\theta^*_t\in L^2(\partial\Omega\times (0,T))$,
$\exists \bar \theta>0: \theta^0(x) \ge \bar\theta$, $\theta^*(x,t) \ge \bar \theta$;}
\item[{\rm (iv)}] $\exists\, 0 < \hat\delta \le \delta < 1/4,\ \exists\, 0<c_\vp< C_\vp$ {\rm such that for all
$p\in \real$ we have that\\$c_\vp \max\{1,|p|\}^{-1-\delta} \le \vp'(p) \le C_\vp \max\{1,|p|\}^{-1-\hat\delta}$;}
\item[{\rm (v)}] {\rm $\exists 0 < c_\mu < C_\mu: c_\mu \le \mu(p) \le C_\mu$ for all $p \in \real$;}
\item[{\rm (vi)}] {\rm $p^0 \in W^{1,2}(\Omega) \cap L^\infty(\Omega)$, $p^*\in L^\infty(\partial\Omega\times (0,T))
\cap L^2(0,T; W^{1,2}(\partial\Omega))$, $p^*_t\in L^2(\partial\Omega\times (0,T))$;}
\item[{\rm (vii)}] {\rm $w^0, \chi^0\in L^\infty(\Omega)$, $\chi^0(x) \in [0,1]$ a.\,e., $\io w^0(x)\dd x = 0$;
\item[{\rm (viii)}] $G\in L^\infty(\Omega\times (0,T))$, $G_t\in L^2(\Omega\times (0,T))$;}
\item[{\rm (ix)}] {\rm $\Omega \subset \real^3$ is a bounded connected set of class $C^{1,1}$,
$\alpha: \partial\Omega \to [0,\infty)$ is Lipschitz continuous, $\omega \in L^\infty(\partial\Omega)$,
$\omega(x) \ge 0$ a.\,e., $\ipo \alpha(x)\dd s(x) >0$, $\ipo \omega(x)\dd s(x) >0$.}
\end{itemize}
\end{hypothesis}

\vspace{2mm}
It is worth noting that it follows from \eqref{lu2} and (vii), using the definition of the functions $G$ and $H$,
that
\begin{equation}
\label{mean}
\io w(x,t)\dd x\,=\,\io w^0(x)\dd x\,=\,0 \quad\forall\, t\in [0,T].
\end{equation}

\vspace{3mm} 
The main result of this paper is the following existence result.

\begin{theorem}\label{t1}
Let Hypothesis \ref{h1} hold true. Then there exists a solution $(p,w,\chi,\theta)$
to the system \eqref{inip}--\eqref{wu5}, \eqref{le7}, with the regularity
\begin{align}
\label{regu1}
&p \in  L^\infty(\Omega\times (0,T)), \quad p_t, \nabla\theta \in L^2(\Omega\times (0,T)), \quad
\nabla p \in L^\infty(0,T;L^2(\Omega)),\\[1mm]
\label{regu2}
&\theta, w_t \in L^{\bar p}(\Omega\times (0,T)), \quad w,\chi_t \in L^\infty(0,T;L^{\bar p}(\Omega))
\,\mbox{ for \,$\bar p < 8+a$},\\[1mm]
\label{regu3}
&\theta_t \in L^2(0,T;W^{-1,q^*}(\Omega)) \mbox{ with \,$q^*>1$\, given by } (\ref{qstar}).
\end{align}
\end{theorem}


The proof of Theorem \ref{t1} will be divided into several steps which each constitutes a
new section in this paper.


\section{Cut-off system}\label{cut}

The strategy for solving Problem \eqref{wu1}--\eqref{wu5} and proving Theorem \ref{t1} is the following: 
we choose a parameter $R>0$ and first solve a cut-off system with the intention to let $R$ tend to infinity.
More precisely, for $R>0$ and $z \in \real$ we denote by $$Q_R(z) = \max\{-R, \min\{z, R\}\}$$ the projection onto
$[-R,R]$, and set $$P_R(z) = z - Q_R(z).$$ We further denote
\be{ce1}
\vp_R(p) = \vp(p) + P_R(p)\,, \quad \Phi_R(p) = \int_0^p \vp_R(\tau) \dd \tau\,,
\quad V_R(p) = p \vp_R(p) - \Phi_R(p)\,,
\ee
and
\be{ce2}
\gamma_R(p,\theta) = \left(1+ (p^2 - R^2)^+\right)\gamma(Q_R(\theta^+)),
\ee

for $p,\theta \in \real$, and replace \eqref{wu1}--\eqref{wu5} by the cut-off system
\ber\label{cu1}
\io\left((\dch (\vp_R(p)+w))_t \eta +\frac{1}{\rho_L}\mu(p)\nabla p\cdot\nabla \eta\right)\dd x
&\!\!=\!\!& \ipo \alpha(x)(p^*-p)\eta \dd s(x),\qquad \\[2mm] \label{cu2}
\nu w_t + \lambda_M w  - p \dch - \beta(Q_R(\theta^+)-\theta_c) &\!\!=\!\!& - G + H_R(t)\quad\mbox{a.~e.},\\[2mm]
 \label{cu4}
 \gamma_R(p,\theta) \chi_t + \partial I(\chi) - (1-\rho^*)(\Phi_R(p)+pw)
&\!\!\ni\!\!& L\left(\frac{Q_R(\theta^+)}{\theta_c}-1\right)\quad\mbox{a.~e.},
\eer

\vspace*{-9mm}
\begin{align}
&\io\left(c_0\theta_t \zeta
+ \kappa(Q_R(\theta^+)) \nabla \theta\cdot\nabla \zeta\right)\dd x \,
- \io\left(\frac{1}{\rho_L}\mu(p)Q_R(|\nabla p|^2) + \gamma_R(p,\theta)\chi_t^2
+  \nu w_t^2\right) \zeta \dd x \nonumber\\ \label{cu5}
&\quad+\io Q_R(\theta^+)\left(\frac{L}{\theta_c}\chi_t + \beta w_t\right) \zeta\dd x
\,\,=\,\, \ipo \omega(x)(\theta^*-\theta)\zeta \dd s(x),
\end{align}
for all test functions $\eta,\zeta \in W^{1,2}(\Omega)$, with
\be{le7a}
H_R(t) = - \frac{1}{|\Omega|}\io (p \dch + \beta(Q_R(\theta^+)-\theta_c) -G)\dd x\,.
\ee

For the system \eqref{cu1}--\eqref{le7a}, we prove the following result.

\begin{proposition}\label{t2}
Let Hypothesis \ref{h1} hold and let $R>0$ be given. Then there exists a solution $(p,w,\chi,\theta)$
to \eqref{cu1}--\eqref{le7a}, \eqref{inip}--\eqref{init} with the regularity
$p,w,\chi,\theta,w_t \in L^q(\Omega; C[0,T])$ for $1\le q < 3$, $p_t,\theta_t \in L^2(\Omega\times (0,T))$,
and $\nabla p, \nabla\theta, \chi_t \in L^\infty(0,T;L^2(\Omega))$.
\end{proposition}

\bpf{Proof of Proposition \ref{t2}}\
Let $$M(p) := \int_0^p \mu(\tau)\dd \tau , \quad K_R(\theta) := \int_0^\theta \kappa(Q_R(\tau^+))\dd \tau,$$
and set $v = M(p)$, $z = K_R(\theta)$. Then \eqref{cu1}--\eqref{cu5} is transformed into the system
\ber\label{ku1}
\io\left((\dch(\vp_R(p)+w))_t \eta + \frac{1}{\rho_L}\nabla v\cdot\nabla \eta\right)\dd x
&\!\!=\!\!& \ipo \alpha(x)(p^*-p)\eta \dd s(x),\qquad \\[2mm] \label{ku2}
\nu w_t + \lambda_M w  - p \dch
- \beta(Q_R(\theta^+)-\theta_c) &\!\!=\!\!& - G + H_R(t)\quad\mbox{a.~e.},\\ \label{ku4}
 \gamma_R(p,\theta) \chi_t + \partial I(\chi) - (1-\rho^*)(\Phi_R(p)+pw)
&\!\!\ni\!\!& L\left(\frac{Q_R(\theta^+)}{\theta_c}-1\right)\quad\mbox{a.~e.}, 
\eer

\vspace*{-9mm}
\begin{align} \nonumber
&\io\left(c_0\theta_t \zeta
+ \nabla z\cdot\nabla \zeta\right)\dd x \,-
 \io\left(\frac{1}{\rho_L}\mu(p)Q_R(|\nabla p|^2) + \gamma_R(p,\theta)\chi_t^2
+  \nu w_t^2\right) \zeta \dd x \\[1mm] \label{ku5}
&\quad+\io Q_R(\theta^+)\left(\frac{L}{\theta_c}\chi_t + \beta w_t\right) \zeta\dd x
\,\,=\,\, \ipo \omega(x)(\theta^*-\theta)\zeta \dd s(x),
\end{align}
which we solve by Galerkin approximations. To this end, let $\{e_k; k=0,1,\dots\}$ denote the complete orthonormal system of eigenfunctions of the problem
\be{eigen}
-\Delta e_k = \lambda_k e_k \ \mbox{ in\ } \Omega \,, \quad \nabla e_k \cdot n = 0
\ \mbox{ on\ } \partial\Omega\,.
\ee
We approximate $v$ and $z$ by the finite sums
\be{ge1}
v\on(x,t) = \sum_{k=0}^n v_k(t) e_k(x)\,, \quad z\on(x,t) = \sum_{k=0}^n z_k(t) e_k(x)\,,
\ee
where $v_k, z_k, w\on, \chi\on$  satisfy the system
\ber \nonumber
&&\hspace{-16mm}\io\left(((\chi\on{+}\rho^*(1{-}\chi\on))(\vp_R(p\on)+w\on))_t e_k
+ \frac{1}{\rho_L}\nabla v\on\cdot\nabla e_k\right)\dd x \\ \label{gu1}
&=& \ipo \alpha(x)(p^*-p\on) e_k \dd s(x), \quad k=0,1, \dots, n,\\[2mm] \nonumber
&&\hspace{-16mm}\nu w\on_t + \lambda_M w\on  - p\on(\chi\on {+} \rho^*(1{-}\chi\on))
- \beta(Q_R((\theta\on)^+)-\theta_c) \\[1mm] \label{gu2}
&=& - G + H\on_R(t)\quad\mbox{a.~e.},\\[2mm] \nonumber
&&\hspace{-16mm} \gamma_R(p\on,\theta\on) \chi\on_t + \partial I(\chi\on)
- (1-\rho^*)(\Phi_R(p\on)+p\on w\on)\\[1mm] \label{gu4}
&\ni& L\left(\frac{Q_R((\theta\on)^+)}{\theta_c}-1\right)\quad\mbox{a.~e.},\\[2mm] \nonumber
&&\hspace{-16mm}\io\left(c_0\theta\on_t e_k + \nabla z\on\cdot\nabla e_k\right)
+  Q_R((\theta\on)^+)\left(\frac{L}{\theta_c}\chi\on_t + \beta w\on_t\right) e_k\dd x \\ \nonumber
&&- \io\left(\frac{1}{\rho_L}\mu(p)Q_R(|\nabla p\on|^2) + \gamma_R(p\on,\theta\on)(\chi\on_t)^2
+  \nu (w\on_t)^2\right) \zeta \dd x \\ \label{gu5}
&=& \ipo \omega(x)(\theta^*-\theta\on)e_k \dd s(x),
\eer
with $p\on := M^{-1}(v\on)$, $\theta\on := K_R^{-1}(z\on)$, and
\be{le7c}
H\on_R(t) := - \frac{1}{|\Omega|}\io (p\on(\chi\on + \rho^*(1{-}\chi\on))
+ \beta(Q_R((\theta\on)^+)-\theta_c) -G)\dd x\,,
\ee
and with the initial conditions
\begin{align} \label{inipn}
v_k(0) &= \io M(p^0(x)) e_k(x)\dd x\,,\\ \label{initn}
z_k(0) &= \io K_R(\theta^0(x)) e_k(x)\dd x\,,\\ \label{iniun}
w\on(x,0) &= w^0(x)\,,\\ \label{inichn}
\chi\on(x,0) &= \chi^0(x)\,.
\end{align}
This is an easy ODE system that admits a unique solution on some interval $[0,T_n) \subset [0,T]$.
Moreover, the solution $w\on$ of \eqref{gu1} enjoys the explicit representation
\ber\nonumber
&& w\on(x,t) \,=\, \expe^{-(\lambda_M/\nu)t} w^0(x)+\frac{1}{\nu}\int_0^t \expe^{(\lambda_M/\nu)(t'-t)}(-G + H\on_R)(x,t')\dd t'\\ \label{n1}
&&\hspace{5mm} +\, \frac{1}{\nu} \int_0^t \expe^{(\lambda_M/\nu)(t'-t)}
\left(p\on(\chi\on {+} \rho^*(1{-}\chi\on)) + \beta(Q_R((\theta\on)^+)-\theta_c)\right)(x,t')\dd t'.\qquad
\eer
Also \eqref{gu4} is of a standard form, namely,
\be{n2}
\chi\on_t + \partial I(\chi\on) \ni F\on,
\ee
with
\be{n3}
F\on = (1-\rho^*) \frac{\Phi_R(p\on)+p\on w\on}{\gamma_R(p\on,\theta\on)}
+ \frac{L(Q_R((\theta\on)^+) -\theta_c)}{\theta_c\gamma_R(p\on,\theta\on)},
\ee
or, equivalently,
\be{n3a}
\chi\on \in [0,1]\,, \quad (F\on - \chi\on_t)(\chi\on - \tilde\chi) \ge 0 \ \mbox{ a.\,e. }\ \forall \tilde\chi \in [0,1]. 
\ee
By virtue of \eqref{n1}--\eqref{n3}, we have for all $(x,t) \in \Omega\times (0,T_n)$ the inequalities
\be{ge3}
\left.
\begin{array}{rcl}
|w\on(x,t)|+|\chi\on_t(x,t)| &\le& C_R \Big(1+\int_0^t|p\on(x,t')|\dd t'
 + \int_0^t\io|p\on(x',t')|\dd x'\dd t'\Big)\\[2mm]
|w\on_t(x,t)| &\le& C_R \Big(1+|p\on(x,t)| + \int_0^t|p\on(x,t')|\dd t'\\
&&+\, \io |p\on(x',t)|\dd x' + \int_0^t\io|p\on(x',t')|\dd x'\dd t'\Big)
\end{array}
\right\}
\ee
where, here and in the following, $C_R>0$ denote constants which possibly depend on $R$ and on the data,
 but  not on $n$.
 
 We now derive some a priori estimates for the solutions to the Galerkin system.
To begin with, we first test \eqref{gu1} by $v_k(t)$ and sum over $k=1, \dots, n$ to obtain 
the identity
(note that $v\on = M(p\on)$,
by definition)
\ber\nonumber
&&\hspace{-12mm} (1-\rho^*) \io \chi\on_t(\vp_R(p\on) + w\on) M(p\on)\dd x
+ \io (\chi\on{+}\rho^*(1{-}\chi\on))\vp_R(p\on)_t M(p\on)\dd x \\ \nonumber
&&+\,\io(\chi\on{+}\rho^*(1{-}\chi\on)) w\on_t M(p\on)\dd x
+ \frac{1}{\rho_L}\io|\nabla v\on|^2 \dd x\\ \label{gn1}
&&+\, \ipo \alpha(x)(p\on-p^*)M(p\on)\dd s(x) = 0.
\eer
We rewrite the first term of \eqref{gn1}, using the identity
\ber\nonumber
\io \chi\on_t(\vp_R(p\on) + w\on) M(p\on)\dd x
&=& \io \chi\on_t(\Phi_R(p\on) + p\on w\on) \frac{M(p\on)}{p\on}\dd x\\ \label{gn2}
&&+\, \io \chi\on_t V_R(p\on)\frac{M(p\on)}{p\on} \dd x.
\eer
From \eqref{gu4} it follows that for a.\,e. $(x,t) \in \Omega\times (0,T)$ we have,
by Young's inequality,
\ber\nonumber
&&\hspace{-12mm}(1-\rho^*)\chi\on_t(\Phi_R(p\on) + p\on w\on) \frac{M(p\on)}{p\on}\\ \nonumber
&=& \frac{M(p\on)}{p\on}\left(
\gamma_R(p\on,\theta\on) \bigl|\chi\on_t\bigr|^2 + \frac{L}{\theta_c}(\theta_c - Q_R((\theta\on)^+))
\chi_t\on\right)\\ \label{gn3}
&\ge& \frac{c_{\mu}}{2} \gamma_R(p\on,\theta\on) \bigl|\chi\on_t\bigr|^2 - C_R\,.
\eer
The second term in \eqref{gn1} can be rewritten as
\ber\nonumber
\io (\chi\on{+}\rho^*(1{-}\chi\on))\vp_R(p\on)_t M(p\on)\dd x &=&
\frac{\dd}{\dd t} \io (\chi\on{+}\rho^*(1{-}\chi\on)) V_{R,M}(p\on)\dd x\\ \label{gn4}
&&-\, (1-\rho^*)\io \chi\on_t V_{R,M}(p\on)\dd x,
\eer
where we denote
\be{gn5}
V_{R,M}(p) = \int_0^p \vp_R'(\tau) M(\tau)\dd \tau\,.
\ee
We see, in particular, that there exist constants $ 0 < c_{R,\mu} < C_{R,\mu}$ such that
$$
c_{R,\mu} p^2 \le V_{R,M}(p) \le C_{R,\mu} p^2
$$
for all $p \in \real$. Combining \eqref{gn2}--\eqref{gn5} with \eqref{gn1}, we obtain that
\ber\nonumber
&&\hspace{-12mm}\frac{\dd}{\dd t} \io (\chi\on{+}\rho^*(1{-}\chi\on)) V_{R,M}(p\on)\dd x
+ \frac{c_{\mu}}{2} \io\gamma_R(p\on,\theta\on) \bigl|\chi\on_t\bigr|^2 \dd x
+ \frac{1}{\rho_L}\io|\nabla v\on|^2 \dd x\\ \nonumber
&&+\, \ipo \alpha(x)(p\on-p^*)M(p\on)\dd s(x)\\ \nonumber
&\le& C_R \,- \io(\chi\on{+}\rho^*(1{-}\chi\on)) w\on_t M(p\on)\dd x\\ \label{gn6}
&&+\, (1-\rho^*)\io \chi\on_t \left(V_{R,M}(p\on) - V_R(p\on)\frac{M(p\on)}{p\on}\right)\dd x.
\eer
By \eqref{ge3} and H\"older's inequality, we have that
\be{gn6a}
\left|\io(\chi\on{+}\rho^*(1{-}\chi\on)) w\on_t M(p\on)(x,t)\dd x\right|
\le C_R \io \left(|p\on|^2(x,t) + \int_0^t |p\on|^2(x,t') \dd t'\right)\dd x,
\ee
and, similarly,
\ber \nonumber
&&\hspace{-12mm}\left|\io \chi\on_t \left(V_{R,M}(p\on) - V_R(p\on)\frac{M(p\on)}{p\on}\right)\dd x\right|
\le C_R \io |\chi\on_t| |p\on|^2 \dd x \\[1mm] \nonumber
&\le& \frac{c_{\mu}}{4} \io\gamma_R(p\on,\theta\on) \bigl|\chi\on_t\bigr|^2 \dd x
+ C_R \io \frac{|p\on|^4}{\gamma_R(p\on,\theta\on)}\dd x \\ [1mm]\label{gn7}
&\le& \frac{c_{\mu}}{4} \io\gamma_R(p\on,\theta\on) \bigl|\chi\on_t\bigr|^2 \dd x
+ C_R \io \left(1+|p\on|^2\right)\dd x.
\eer
Let $[0,T_n)$ denote the maximal interval of existence of our solution.
Using \eqref{gn6}--\eqref{gn7}, and Gronwall's lemma, we thus can infer that
\be{gn8}
\supess_{t\in (0,T_n)}\io |p\on|^2(x,t)\dd x +\int_0^{T_n}\io |\nabla p\on|^2\dd x\dd t
+ \int_0^{T_n}\ipo \alpha(x)|p\on|^2 \dd s(x)\dd t \le C_R\,.
\ee
In particular, the Galerkin solution exists globally, and for every $n \in \nat$ we have $T_n=T$.

In what follows, we denote by $|\cdot|_p$ the norm in $L^p(\Omega)$, by $\|\cdot\|_p$ the norm
in $L^p(\Omega\times (0,T))$, by $\|\cdot\|_{\partial\Omega, p}$ the norm
in $L^p(\partial\Omega\times (0,T))$, and by $\|\cdot\|_{W^{\ell,p}(\Omega)}$ the norm in $W^{\ell,p}(\Omega)$
for $\ell \in \nat$ and $1 \le p \le \infty$.

Let us recall the Gagliardo-Nirenberg inequality
\be{gn}
|u|_q \le C\left(|u|_s + |u|_s^{1-\rho}|\nabla u|_p^\rho\right),
\ee
with
$$
\rho = \frac{\frac 1s - \frac1q}{\frac 1s + \frac 1N - \frac 1p}\,\,,
$$
which is valid for all $1\le s<q$, $1/q > (1/p) - (1/N)$, every bounded open set $\Omega \subset \real^N$
with Lipschitzian boundary, and every function $u \in W^{1,p}(\Omega)$.
For $t\in (0,T)$, $N=3$, $s = p = 2$, and $q=4$, we have, in particular,
$$
|p\on(t)|_4 \le C\left(|p\on(t)|_2 + |p\on(t)|_2^{1/4} |\nabla p\on(t)|_2^{3/4}\right).
$$
Hence, by \eqref{gn8},
\be{gn9}
\int_0^T \left(\io |p\on(x,t)|^4 \dd x\right)^{2/3}\dd t \le C_R\,, 
\ee
independently of $n$.

Next, we test \eqref{gu1} by $\dot v_k(t)$ and sum over $k=0,1,\dots n$ to obtain the identity
\ber\nonumber
&&\hspace{-16mm}\frac{\dd}{\dd t} \left(\frac1{2\rho_L}\io |\nabla v\on|^2\dd x
+ \ipo \alpha(x)(\hat M(v\on)-p^* v\on)\dd s(x)\right)\\ \nonumber
&& +\io\left(((\chi\on{+}\rho^*(1{-}\chi\on))(\vp_R(p\on)+w\on))_t v_t\on\right)\dd x\\ \label{gd0}
&=& - \ipo \alpha(x) p^*_t v\on \dd s(x),
\eer
where $\hat M' = M^{-1}$.
We have the pointwise lower bound
$$
\vp_R(p\on)_t v\on_t \ge C_R |v\on_t|^2,
$$
and thus, by \eqref{gd0} and H\"older's inequality, we have for all $t \in [0,T]$ that 
\ber\nonumber
&&\hspace{-16mm}\int_0^t\io |v\on_t|^2\dd x\dd t' + \left(\io |\nabla v\on|^2(x,t) \dd x
+ \ipo \alpha(x)|v\on|^2(x,t)\dd s(x)\right) \\ \label{gd1}
&\le& C_R \left(1+ \int_0^t\io (|w\on_t|^2 + |\chi\on_t|^2 |w\on|^2)\dd x\dd t'
+\int_0^t\ipo \alpha(x)|v\on|^2\dd s(x)\dd t'\right).\quad
\eer
A bound for $\|w\on_t\|_2^2$ follows from \eqref{ge3} and \eqref{gn8}. Moreover, owing to 
\eqref{ge3}, we have for $t\in (0,T)$
that
$$
\io |\chi\on_t|^2 |w\on|^2 (x,t)\dd x \le C_R\left(1 + \io \left(\int_0^t |p\on(x,t')|\dd t'\right)^4\dd x\right).
$$
We use the Minkowski inequality in the form
$$
\left(\io \left(\int_0^t |p\on(x,t')|\dd t'\right)^4\dd x\right)^{1/4} \le
\int_0^t \left(\io |p\on(x,t')|^4\dd x\right)^{1/4} \dd t'
$$
to check that
$$
\io |\chi\on_t|^2 |w\on|^2 (x,t)\dd x
\le C_R\left(1 + \left(\int_0^t \left(\io |p\on(x,t')|^4\dd x\right)^{1/4} \dd t'\right)^4\right) \le C_R,
$$
by virtue of \eqref{gn9}. Then \eqref{gd1} and the Gronwall argument imply that
\be{ge2}
\|v\on(t)\|_{W^{1,2}(\Omega)}^2 + \int_0^t |v\on_t(t')|_2^2\dd t' \le C_R,
\ee
whence also
\be{ge2a}
\|p\on(t)\|_{W^{1,2}(\Omega)}^2 + \int_0^t |p\on_t(t')|_2^2\dd t' \le C_R
\ee
for $t \in (0,T)$. 

We continue by testing \eqref{gu5} by $\dot z_k(t)$ and summing over $k=0,1,\dots n$.
Note that, thanks to \eqref{n2}--\eqref{n3}, \eqref{ge3}, and \eqref{gn8}, we have that
$$
\gamma_R(p\on,\theta\on)(\chi\on_t(x,t))^2 +  \nu (w\on_t(x,t))^2 \le 
C_R \left(1+|p\on(x,t)| + \int_0^t|p\on(x,t')|\dd t'\right)^3
$$
for a.\,e. $(x,t) \in \Omega \times (0,T)$. This yields the inequality
\bers
&&\hspace{-16mm}\frac{\dd}{\dd t} \left(\frac12\io |\nabla z\on|^2\dd x
+ \ipo \omega(x)(\hat K_R(v\on)-\theta^* v\on)\dd s(x)\right) +
\frac12 \io c_0\theta\on_t z\on_t \dd x\\
&\le&  \ipo \omega(x) |\theta^*_t| |z\on| \dd s(x)
+ C_R\io \left(1+|p\on(x,t)| + \int_0^t|p\on(x,t')|\dd t'\right)^6\dd x,
\eers
where $\hat K_R' = K_R^{-1}$. Using \eqref{ge2a} and the Sobolev embedding theorem, we obtain, as before, that
\be{ge2b}
\|\theta\on(t)\|_{W^{1,2}(\Omega)}^2 + \int_0^t |\theta\on_t(t')|_2^2\dd t' \le C_R
\ee
for $t \in (0,T)$. 
Hence, there exist a subsequence of $\{(p\on, \theta\on): n \in \nat\}$, which is again indexed by $n$,
and functions $p, \theta$, such that
\bers
p\on_t \to p_t, \ \theta\on_t \to \theta_t, && \mbox{weakly in } \ L^2(\Omega\times (0,T)),\\
\nabla p\on \to \nabla p, \ \nabla \theta\on \to \nabla \theta, && \mbox{weakly-star in } \ L^\infty(0,T;L^2(\Omega)),\\
p\on \to p, \ \theta\on \to \theta, && \mbox{strongly in } \ L^q(\Omega; C[0,T]) \ \for \ 1\le q < 3,\\
\eers
where we used the compact embedding
$W^{1,2}(\Omega\times (0,T))\hookrightarrow\hookrightarrow L^q(\Omega;C[0,T])$ for $1\le q<3$, see \cite{bin}.

We now check that the sequences $\{w\on\}, \{\chi\on\}, \{w\on_t\}, \{\chi\on_t\}$ converge strongly in appropriate function spaces and that the limit functions
satisfy the system \eqref{cu1}--\eqref{le7a}. Passing again to a subsequence if necessary, we may fix a set
$\Omega' \subset \Omega$ with meas$(\Omega\setminus\Omega') = 0$ such that
\be{n4}
\lim_{n\to \infty} \sup_{t\in [0,T]}|p\on(x,t) - p(x,t)| = 0\,, \quad \lim_{n\to \infty} \sup_{t\in [0,T]}|\theta\on(x,t) - \theta(x,t)| = 0,
\quad \forall\, x \in \Omega',
\ee
and such that the functions $t \mapsto p(x,t)$ and $t \mapsto \theta(x,t)$ belong to $C[0,T]$ for all $x \in \Omega'$. In particular, we can define the
real numbers
$$
\widetilde p(x) := \sup_{t\in [0,T]}|p(x,t)|\,, \quad \widetilde \theta(x): = \sup_{t\in [0,T]}|\theta(x,t)|,
 \ \for x \in \Omega'\,.
$$
Let $x \in \Omega'$ be arbitrarily fixed now. Then there is some $n_0(x) \in \nat$ such that for $n > n_0(x)$ we have
$|p\on(x,t)| \le 2\widetilde p(x)$ and $|\theta\on(x,t)| \le 2\widetilde \theta(x)$, for all $t \in [0,T]$ and $x \in \Omega'$.

For $n,m \in \nat$, $n,m > n_0(x)$, we have by \eqref{n1} for $t \in [0,T]$ and $x \in \Omega'$ that
\be{n5}
|w\on(x,t) - w\om(x,t)| \le C_R(1+\widetilde p(x))\int_0^t (|\chi\on - \chi\om| + |p\on - p\om| + |\theta\on - \theta\om|)(x,t')\dd t'.
\ee
Hence, with the notation of \eqref{n3},
\ber\nonumber
&&\hspace{-16mm}\int_0^t|F\on(x,t') - F\om(x,t')|\dd t'\\ \label{n6}
&\le& C_R(1+\widetilde p(x))^2\int_0^t (|\chi\on - \chi\om| + |p\on - p\om| + |\theta\on - \theta\om|)(x,t')\dd t'.
\eer
The well-known $L^1$-Lipschitz continuity result for variational inequalities (see, e.\,g., \cite[Theorem 1.12]{cmuc})
tells us that 
\be{n7}
\int_0^t |\chi\on_t - \chi\om_t|(x,t')\dd t' \le 2 \int_0^t|F\on(x,t') - F\om(x,t')|\dd t'.
\ee
Since $\{p\on(x,\cdot)\}$ and $\{\theta\on(x,\cdot)\}$ converge uniformly for each $x \in \Omega'$, we may apply the Gronwall argument
to conclude that $\{\chi\on(x,\cdot)\}$, $\{w\on(x,\cdot)\}$, $\{w_t\on(x,\cdot)\}$ are Cauchy sequences in $C[0,T]$
and that $\{\chi_t\on(x,\cdot)\}$ is a Cauchy sequence in $W^{1,1}(0,T)$, for every $x \in \Omega'$. Hence, there exist functions
$\chi, w : \Omega'\times (0,T)$ such that, as $n\to\infty$,
\be{n8}
\sup_{t \in [0,T]}|w\on(x,t) - w(x,t)| \to 0\,, \ \ \sup_{t \in [0,T]}|\chi\on(x,t) - \chi(x,t)| \to 0\,, \ \
\sup_{t \in [0,T]}|w_t\on(x,t) - w_t(x,t)| \to 0
\ee
and
\be{n9}
\int_0^T|\chi\on_t(x,t) - \chi_t(x,t)| \dd t \to 0,
\ee
for all $x \in \Omega'$. Since $|\chi\on|, |w\on|, |\chi\on_t|, |w\on_t|$ admit a pointwise upper
bound \eqref{ge3} in terms of convergent sequences in $L^q(\Omega; C[0,T])$ for $1\le q < 3$, we can use
the Lebesgue Dominated Convergence Theorem to conclude that
\be{n10}
w\on \to w\,, \quad w_t\on \to w_t\,, \quad \chi\on \to \chi, \quad \mbox{ strongly in } \ L^q(\Omega; C[0,T])\,,
\ee 
and 
\be{n11}
\chi\on_t \to \chi_t \quad \mbox{ strongly in } \ L^1(\Omega\times (0,T)).
\ee 
Moreover, from H\"older's inequality we obtain that
\bers
&&\hspace{-12mm}\int_0^T\io |\chi\on_t - \chi_t|^2\dd x \dd t = \int_0^T\io |\chi\on_t - \chi_t|^{1/3} |\chi\on_t - \chi_t|^{5/3}\dd x \dd t\\
&& \le\, \left(\int_0^T\io |\chi\on_t - \chi_t|\dd x \dd t\right)^{1/3} \left(\int_0^T\io |\chi\on_t - \chi_t|^{5/2}\dd x \dd t\right)^{2/3}
\to 0,
\eers
by virtue of \eqref{n11} and \eqref{ge3}. We can therefore pass to the limit
in \eqref{gu1}--\eqref{le7c}, where \eqref{gu4} is interpreted as the variational inequality \eqref{n3a},
and check that its limit is the desired solution to \eqref{cu1}--\eqref{le7a}.
\epf


\section{Estimates independent of $R$}\label{apr}
In this section, we derive estimates for the solutions to \eqref{cu1}--\eqref{cu5} which are independent of the
cut-off parameter $R$. In the entire section, we denote by $C$ positive constants which
may depend on the data of the problem, but not on $R$.

\subsection{Positivity of temperature}\label{pos}

For every nonnegative test function $\zeta \in W^{1,2}(\Omega)$ we have, by virtue of \eqref{cu5},
and using the fact that $\gamma_R(p,\theta) \ge c_\gamma>0$,
\be{pe25}
\io\left(c_0\theta_t \zeta + \kappa(Q_R(\theta^+)) \nabla \theta\cdot\nabla \zeta\right)\dd x
+ \ipo \omega(x)(\theta-\theta^*) \zeta \dd s(x) \ge - C\io (Q_R(\theta^+))^2 \zeta \dd x
\ee
with a constant $C$ depending only on the constants $L, \theta_c, \beta, \nu, c_\gamma$.
Let $\psi$ be the solution of the equation
\be{epsi}
c_0 \dot\psi(t) + C\psi^2(t) = 0\,, \quad \psi(0) = \bar\theta\,.
\ee
Then
\be{epsi2}
\psi(t) = \frac{\bar\theta c_0}{c_0 + \bar\theta C t}\,\,,
\ee
and we have
\ber \nonumber
&&\hspace{-16mm}\io\left(c_0(\psi- \theta)_t \zeta
+ \kappa(Q_R(\theta^+)) \nabla(\psi- \theta)\cdot\nabla \zeta\right)\dd x\,
- \ipo \omega(x)(\theta-\theta^*)\zeta \dd s(x)\\ \label{pe25b}
&\le&  C\io ((Q_R(\theta^+))^2 - \psi^2) \zeta \dd x
\eer
for every nonnegative test function $\zeta \in W^{1,2}(\Omega)$. In particular, for
$\zeta(x,t) = (\psi(t) - \theta(x,t))^+$, we obtain that
\be{pe25c}
\frac{\dd}{\dd t}\frac{c_0}{2}\io((\psi- \theta)^+)^2\dd x 
+ \ipo \omega(x)(\theta^*-\theta)(\psi- \theta)^+ \dd s(x)
\le  C\io ((Q_R(\theta^+))^2 - \psi^2)(\psi- \theta)^+ \dd x\,.
\ee
From Hypothesis 2.1\,(iii), we obtain for all values of $x$ and $t$ that
$$
(\theta^*-\theta)(\psi- \theta)^+ \ge 0\,,
\quad ((Q_R(\theta^+))^2 - \psi^2)(\psi- \theta)^+
= (Q_R(\theta^+) - \psi)(Q_R(\theta^+) + \psi)(\psi- \theta)^+ \le 0\,,
$$
and from 
\be{pe25d}
\frac{\dd}{\dd t}\frac{c_0}{2}\io((\psi- \theta)^+)^2\dd x \le 0\,, \quad (\psi- \theta)^+(x,0) = 0,
\ee
we conclude that, independently of $R>0$,
\be{pos0}
\theta(x,t) \ge \psi(t) \ge \frac{\bar\theta c_0}{c_0 + \bar\theta C T} > 0
\quad \mbox{ for all }\ x \ \mbox{ and } \ t.
\ee


\subsection{Energy estimate}\label{enes}

We test \eqref{cu1} by $\eta=p$, \eqref{cu5} by $\zeta = 1$, and sum up. With the notation \eqref{ce1}, we use
the identities
\ber\label{en5}
\io(\dch\vp_R(p))_t p \dd x &=& \frac{\dd}{\dd t}\io \dch V_R(p)\dd x
+ \io(1-\rho^*)\Phi_R(p)\chi_t \dd x\,,\qquad\\ \label{en5a}
\io(\dch w)_t p \dd x &=& \io \dch w_t p \dd x
+ \io(1-\rho^*)wp\chi_t \dd x\,,\qquad\\ \label{en6}
\chi_t \left(\gamma_R(p,\theta)\chi_t - L\frac{Q_R(\theta)}{\theta_c}\right)
&=& -L\chi_t + (1-\rho^*)(\Phi_R(p)+pw)\chi_t,
\eer
which follow from \eqref{cu4}, and we obtain that
\ber \nonumber
&&\frac{\dd}{\dd t}\io \left(c_0\theta + L\chi + \dch V_R(p)+ \frac{\lambda_M}{2} w^2
+ \beta\theta_c w\right)\dd x \\ \label{en1}
&& \qquad + \ipo \left(\omega(x)(\theta-\theta^*) + \alpha(x)(p-p^*)p\right)\dd s(x) \le
\io w_t(H_R(t) - G)\dd x.
\eer
Note that $V_R'(p) = p \vp_R'(p)$, $V_R(0) = 0$, so that $V_R(p) > 0$ for all $p \ne 0$. Furthermore,
$$
\io w_t H_R(t)\dd x = H_R(t) \io w_t \dd x = 0,
$$
so that \eqref{en1} can be written as
\ber \nonumber
&&\frac{\dd}{\dd t}\io \left(c_0\theta + L\chi + \dch V_R(p)+ \frac{\lambda_M}{2} w^2
+ (\beta\theta_c+G) w\right )\dd x \\ \label{en1a}
&& \qquad + \ipo \left(\omega(x)(\theta-\theta^*) + \alpha(x)(p-p^*)p\right)\dd s(x) \le  
\io G_t w\dd x.
\eer
By Gronwall's argument and Hypothesis \ref{h1}, we thus have
\be{en0}
\supess_{t\in (0,T)} \io (\theta + V_R(p) + w^2)\dd x
+ \int_0^T\ipo(\omega(x)\theta + \alpha(x) p^2)\dd s(x)\dd t \le C\,.
\ee


\subsection{The Dafermos estimate}\label{dafe}

We denote $\hate = Q_R(\theta) = Q_R(\theta^+)$ and rewrite \eqref{cu5} in the form
\begin{align}
&\io\left(c_0\theta_t \zeta
+ \kappa(\hate) \nabla \theta\cdot\nabla \zeta\right)\dd x\,
- \io\left(\frac{1}{\rho_L}\mu(p)Q_R(|\nabla p|^2) + \gamma_R(p,\theta)\chi_t^2
+  \nu w_t^2\right) \zeta \dd x \nonumber\\ \label{cu5a}
&+\io \hate\left(\frac{L}{\theta_c}\chi_t + \beta w_t\right) \zeta\dd x
\,=\, \ipo \omega(x)(\theta^*-\theta)\zeta \dd s(x),
\end{align}
for every $\zeta \in W^{1,2}(\Omega)$. We test \eqref{cu5a} by $\zeta = -\hate^{-a}$. This yields
the identity
\ber \nonumber
&&\hspace{-12mm}\io\frac{a\kappa(\hate)}{\hate^{1+a}} |\nabla \hate|^2\dd x +
 \io\hate^{-a}\left(\frac{1}{\rho_L}\mu(p)Q_R(|\nabla p|^2) + \gamma_R(p,\theta)\chi_t^2
+  \nu w_t^2\right) \dd x \\ \label{cu5b}
&=&\io \hate^{1-a}\left(\frac{L}{\theta_c}\chi_t + \beta w_t\right) \dd x
+\ipo \omega(x)(\theta-\theta^*)\hate^{-a} \dd s(x) + \frac{c_0}{1-a}\frac{\dd}{\dd t}\io\hate^{1-a}\dd x.\qquad
\eer
By Hypothesis \ref{h1}\,(ii), we have $\frac{\kappa(\hate)}{\hate^{1+a}} \ge c_\kappa$. Furthermore, H\"older's and Young's inequalities give the estimate
\be{de0}
\io \hate^{1-a}\left(|\chi_t| + |w_t|\right) \dd x \le \frac C\tau \io \hate^{2-a}\dd x + \tau \io
\hate^{-a}\left(\chi_t^2 +  w_t^2\right) \dd x
\ee
for every $\tau > 0$. This and \eqref{en0} yield the estimate
\be{de1}
\int_0^T\io|\nabla\hate(t)|^2\dd x\dd t \le C\left(1 +  \int_0^T\io \hate^{2-a}\dd x\dd t\right)\,.
\ee
{}From the Gagliardo-Nirenberg inequality \eqref{gn} with $s=1$, $p=2$, and $N=3$, we obtain that
\be{de1a}
|\hate(t)|_q \le C\left(1+ |\nabla\hate(t)|_2^\rho\right),
\ee
with $\rho = (6/5(1 - (1/q))$, where we used \eqref{en0} once more. In particular, for every $q \le 8/3$, we have by \eqref{de1} and \eqref{de1a} that
\be{de2a}
\int_0^T|\hate(t)|_q^{5q/3(q-1)}\dd t \le C\left(1+ \int_0^T|\nabla\hate(t)|_2^2\dd t\right)
\le C\left(1 +  \int_0^T|\hate|_{2-a}^{2-a}\dd x\dd t\right)\,.
\ee
Moreover, using \eqref{de2} first for $q = 2-a$ and then for $q=8/3$, we obtain that
\be{de2}
\int_0^T|\hate(t)|_{8/3}^{8/3}\dd t + \int_0^T|\nabla\hate(t)|_2^2\dd t \le C,
\ee
independently of $R$.


\subsection{Estimates for the capillary pressure}\label{capi}

We choose an even function $b:\real \to (0,\infty)$ such that the functions $b$ and $p \mapsto p b(p)$ are
Lipschitz continuous and such that $b'(p) \ge 0$ for $p > 0$. Then, owing to \eqref{ge2a}, $\eta = p b(p)$ is an admissible
test function in \eqref{cu1}, and the term under the time derivative has the form
\ber \nonumber
&&\io (\dch (\vp_R(p)+w))_t p b(p) \dd x \,= \io \dch \vp_R(p)_t p b(p) \dd x\\ \label{bbe1}
&& \quad +\,\io \dch w_t  p b(p) \dd x
+ (1-\rho^*)\io \chi_t \left(\varphi_R(p) + w\right)p b(p) \dd x.
\eer
For $p \in \real$, we put
\begin{align*}
V_b(p) &\,:=\, \int_0^p \vp'(\tau)\, \tau\, b(\tau) \dd \tau,\quad
\hat P_{R,b}(p) \,:=\, \int_0^p P'_R(\tau)\, \tau\, b(\tau) \dd\tau,\\
\Psi_{R,b}(p) &\,:=\, \vp_R(p) p b(p) - \hat P_{R,b}(p) - V_b(p) - \Phi_R(p) b(p) \,=\,
\int_0^p V_R(\tau) \,\tau\, b'(\tau) \dd\tau\,.
\end{align*}
Then $V_b(p) > 0$ and $\Psi_{R,b}(p) \ge 0$ for all $p\ne 0$, and \eqref{bbe1} can be rewritten as
\ber \nonumber
&&\io (\dch (\vp_R(p)+w))_t p b(p) \dd x = \frac{\dd}{\dd t} \io \dch (\hat P_{R,b}(p) + V_b(p)) \dd x\\ \label{bbe2}
&& \quad +\,\io \dch w_t  p b(p) \dd x
+ (1-\rho^*)\io \chi_t \left((\Phi_R(p) + wp)b(p) + \Psi_{R,b}(p)\right)\dd x.\qquad
\eer
Owing to \eqref{en6}, we have, with the notation from Subsection \ref{dafe}, that
$$
(1-\rho^*)\chi_t (\Phi_R(p) + wp)\, =\, \gamma_R(p,\theta) \chi_t^2 \,+\, \frac{L}{\theta_c} (\theta_c - \hate)
\chi_t\,
\ge\, \frac12 \gamma_R(p,\theta) \chi_t^2 - C (1+ \hate)
$$
with a constant $C>0$ independent of $R$. Similarly,
$$
\left|\io \chi_t \Psi_{R,b}(p)\dd x \right| \le \frac14 \io \gamma_R(p,\theta) \chi_t^2 b(p) \dd x
+ C \io \frac{\Psi_{R,b}^2(p)}{b(p)\gamma_R(p,\theta)}\dd x.
$$
We have, by definition, that $\Psi_{R,b}(p) \le V_R(p) b(p)$, hence
$$
\frac{\Psi_{R,b}^2(p)}{b(p)\gamma_R(p,\theta)} \,\le\, C\frac{V_{R}^2(p) b(p)}{\gamma_R(p,\theta)}\,
\le\, C b(p)\frac{(V(p) + \frac12(p^2 - R^2)^+)^2}{1 + (p^2 - R^2)^+} \,\le\, C p^2 b(p)
$$
independently of $R$. We conclude that
\ber \nonumber
&&\io (\dch (\vp_R(p)+w))_t p b(p) \dd x
\ge \frac{\dd}{\dd t} \io \dch (\hat P_{R,b}(p) + V_b(p)) \dd x\\ \label{bbe3}
&& \quad +\,\frac14 \io \gamma_R(p,\theta) \chi_t^2 b(p) \dd x
- C \io \left(1+ |w_t| + |p| + \hate\right) |p| b(p) \dd x.\qquad
\eer
{}From \eqref{cu1}, with $\eta = p b(p)$, we thus obtain, in particular, that
\ber \nonumber 
&&\hspace{-10mm}
\frac{\dd}{\dd t} \io \dch (\hat P_{R,b}(p)+V_b(p)) \dd x + \io \mu(p)(p b'(p) + b(p))|\nabla p|^2 \dd x \\ \label{bbe4}
&& +\, \ipo \alpha(x) (p-p^*) p b(p)\dd s(x)
\le  C \io \left(1{+}|w_t|{+}|p|{+}\hate\right) |p| b(p) \dd x\,,
\eer
with a constant $C>0$ which is independent of both $b$ and $R$.
To estimate the right-hand side of \eqref{bbe4}, we first notice that $|H(t)| \le C(1 + \io |p|\dd x)$,
and from \eqref{cu2} and Hypothesis \ref{h1}\,(viii) we obtain the pointwise bounds
\ber \label{c3}
|w(x,t)| &\le& C\left(1 + \int_0^t(|p(x,{t'})| + \hate(x,{t'}))\dd{t'}
+ \int_0^t\io|p(x',{t'})|\dd x'\dd{t'}\right),\\ \label{c3a}
|w_t(x,t)| &\le& |w(x,t)|+ C\left(1 + |p(x,t)| + \hate(x,t)+ \io|p(x',t)|\dd x' \right)\,. 
\eer
In particular, for $b(p) \equiv 1$ we have $\Psi_{R,b}(p) = 0$ and $V_b=V,$ and it follows from \eqref{bbe4}--\eqref{c3}
that
\ber\label{esti} \nonumber
&&\hspace{-16mm} \io (\dch V(p))(x,t)\dd x + \int_0^t\io \mu(p)|\nabla p|^2(x,t')\dd x\dd{t'}
+ \int_0^t\ipo \alpha(x) |p|^{2}(x,t')\dd s(x)\dd{t'}\\ \label{c4}
&\le&  C \left(1+ \int_0^t\io \left( \hate |p| + |p|^2\right)\dd x\dd{t'}\right).
\eer
We have, by Hypothesis \ref{h1}\,(iv), that $V(p) \ge c_{\vp}(|p|^{1-\delta} - \delta)/(1-\delta)$.
The energy estimate \eqref{en0} then yields that
\be{c4a}
\io |p|^{1-\delta}(x,t)\dd x \le C\,.
\ee
Moreover, by \eqref{de2}, $\hate$ is bounded in $L^{8/3}(\Omega\times (0,T))$.
We thus obtain from \eqref{c4} that
\ber \nonumber
&&\hspace{-16mm} \io (\dch V(p))(x,t)\dd x + \int_0^t\io \mu(p)|\nabla p|^2\dd x\dd{t'}
+ \int_0^t\ipo \alpha(x) |p|^{2}\dd s(x)\dd{t'}\\ \label{c5}
&\le&  C \left(1+ \int_0^t\io |p|^2\dd x\dd{t'}\right)\,.
\eer
Furthermore, by Hypothesis \ref{h1}\,(ix), $\Omega$ is connected, and $\int_{\partial\Omega}\alpha(x)\dd s(x)>0$. This implies that there exists a constant $C_\Omega>0$, which depends only on $\Omega$, such that, a.~e. in $(0,T)$,
\be{c6}
C_\Omega\,\|p\|_{W^{1,2}(\Omega)}^2\,\le\,\ipo \alpha(x) |p|^{2}\dd s(x) + \io \mu(p)|\nabla p|^2\dd x\,.
\ee
Moreover, we infer from H\"older's inequality that
\be{c6a}
\io |p|^2\dd x = \io |p|^{(1-\delta)/2}|p|^{(3+\delta)/2}\dd x \le \left(\io |p|^{1-\delta}\dd x\right)^{1/2}
\left(\io |p|^{3+\delta}\dd x\right)^{1/2}\,.
\ee
Hence, by \eqref{c4a},
\begin{align}
\label{c6b}
&\int_0^t\io |p|^2\dd x\dd{t'} \le  C\int_0^t\left(\io |p|^{3+\delta}\dd x\right)^{1/2}\dd{t'}
\le  C\left(\int_0^t\left(\io |p|^{3+\delta}\dd x\right)^{2/(3+\delta)}\dd{t'} \right)^{(3+\delta)/4}
\nonumber
\\
&\quad\le\,C\left(\int_0^t |p(t')|_{3+\delta}^2\dd t'\right)^{(3+\delta)/4}\,.
\end{align}
Since $\delta<1$, we have the embedding inequality
$$
|p(t)|_{3+\delta}^2 \le C\,\|p(t)\|_{W^{1,2}(\Omega)}^2,
$$
so that from \eqref{c6b} it follows that
\be{c7}
\int_0^t\io|p|^2\dd x\dd {t'} \,\le\, C\left(\int_0^t\|p({t'})\|_{W^{1,2}(\Omega)}^2\dd {t'}\right)^{(3+\delta)/4}.
\ee
Employing Young's inequality, we therefore conclude from \eqref{esti} and \eqref{c6} that
\be{c7a}
\|p\|_{L^2(0,T;W^{1,2}(\Omega))} \le C.
\ee
Moreover, for $(x,t) \in \Omega\times (0,T)$,
$q\ge 1$ and $s > 1$, we have
$$
|p(x,t)|^q = |p(x,t)|^{(1-\delta)/s} |p(x,t)|^{q-(1-\delta)/s},
$$
whence, by H\"older's inequality with $s' = s/(s-1)$,
$$
|p(t)|_q^q = \left(\io|p(x,t)|^{1-\delta}\dd x\right)^{1/s}  \left(\io|p(x,t)|^{(q-(1-\delta)/s)s'}\dd x\right)^{1/s'}.
$$
We thus obtain from \eqref{c4a} and \eqref{c7a} that
\be{c8}
\int_0^T |p(t)|_q^q \dd t \le C,
\ee
provided that $(q-(1-\delta)/s)s' \le 6$ and $s' \ge 3$. In other words,
$$
q \le \frac{1-\delta}{s} + \frac{6}{s'} = \frac{5+\delta}{s'} + 1-\delta \le \frac{5+\delta}{3} + 1-\delta,
$$
and the maximal admissible value for $q$ in \eqref{c8} is given by
\be{c9}
q = \frac{8 - 2\delta}{3}\,.
\ee
Let now the function $b$ in \eqref{bbe4} be arbitrary. For $p\in \real$, we put $\hat b(p): = \int_0^p \tau\,b(\tau) \dd \tau$.
Then $\hat b$ is convex, and we have the inequality
\be{hatb}
\hat b(p) - \hat b(p^*) \le (p-p^*) \hat b'(p) = (p-p^*) p b(p).
\ee
{}From \eqref{bbe4}, \eqref{de2}, \eqref{c3}, \eqref{c8}, \eqref{c9}, and \eqref{hatb} it follows that
there exists a function $h\in L^q(\Omega\times (0,T))$ such that
$$
\|h\|_q \le C,
$$
with a constant $C>0$ independent of $b$ and $R$, as well as
\ber \nonumber 
&&\hspace{-10mm}
\frac{\dd}{\dd t} \io \dch(\hat P_{R,b}(p)+ V_b(p)) \dd x + c_{\mu} \io b(p)|\nabla p|^2 \dd x 
+\ipo\alpha(x)\hat b(p)\dd s(x)\\ \label{bbe5}
&\le&  \ipo\alpha(x)\hat b(p^*)\dd s(x) + \io h |p| b(p) \dd x.
\eer

Integration of \eqref{bbe5} in time, using the fact that $\chi+\rho^*(1-\chi)\ge\rho^*>0$, yields the estimate
\begin{align} \nonumber 
&
\io V_b(p)(x,t) \dd x + \int_0^T\io b(p)|\nabla p|^2 \dd x \dd {t'}
+\int_0^T\ipo\alpha(x)\hat b(p)\dd s(x)\dd {t'}\\ \label{bbe6}
&\le\, C\left(\io V_b(p)(x,0) \dd x + \int_0^T\ipo\alpha(x)\hat b(p^*)\dd s(x)\dd{t'}
+ \left(\int_0^T\io (|p| b(p))^{q'} \dd x \dd {t'}\right)^{1/q'}\right),
\end{align}
for all $t \in [0,T]$, with $q' = \frac{q}{q-1} = \frac{8-2\delta}{5-2\delta}$.

Now let $k>0$ be given, and let $\{b_n\}_{n \in \nat}$ be a sequence of even, smooth, bounded
functions which are nondecreasing in $(0,\infty)$ and
such that $b_n(p) \nearrow |p|^{2k}$ locally uniformly in $\real$.
Then $(|p| b_n(p))^{q'} \nearrow |p|^{(1+2k)q'}$
locally uniformly. From \eqref{c8} we know that $p\in L^q(\Omega\times (0,T))$, where $q$ is 
given by \eqref{c9}.
Hence, the integral on the right-hand side of \eqref{bbe6} is meaningful if $(1+2k)q'\le q$, that is,
if $3k \le 1 - \delta$. In particular, thanks to Hypothesis \ref{h1}\,(iv), $k=\delta$ is an admissible choice.

We continue by induction. To this end, assume that
\be{bbe7}
\int_0^T \io |p|^{(1+2k)q'} \dd x \dd t =: J_k < \infty
\ee
holds true for some $k \ge \delta$. Using the denotations
$$V_{b_n}(p):=\int_0^p\varphi'(\tau)\,\tau\,b_n(\tau)\dd\tau,\quad \hat b_n(p):=\int_0^p\tau\,b_n(\tau)\dd\tau \quad
\mbox{for $\,n\in\nat$},
$$
we can estimate the terms occurring on the right-hand side of \eqref{bbe6} for $n\in\nat$ as follows:
\bers
\io V_{b_n}(p)(x,0) \dd x &\le& C|p^0|_\infty^{2k + 1 - \hat\delta}\,,\\
\int_0^T\ipo\alpha(x)\hat b_n(p^*)\dd s(x)\dd{t'} &\le& C\|p^*\|_{\partial\Omega,\infty}^{2k+2}\,.
\eers
Put $E = \max\{1, |p^0|_\infty, \|p^*\|_{\partial\Omega,\infty}\}$. Then \eqref{bbe6}
can for $n\in\nat$ be rewritten as
\be{bbe8}
\io V_{b_n}(p)(x,t) \dd x + \int_0^T\!\!\!\io b_n(p)|\nabla p|^2 \dd x \dd {t'}
+\int_0^T\!\!\!\ipo\!\!\alpha(x)\hat b_n(p)\dd s(x)\dd {t'}
\le C\max\{E^{2k+2}, J_k^{1/q'}\},
\ee
independently of $k$, $R$, and $n$. By virtue of Fatou's lemma, we can take the limit as  $n \to \infty$ to obtain
that \eqref{bbe8} holds true for $b(p) = |p|^{2k}$. Using the estimate
\be{bbe9}
V_b(p) \ge \frac{c_{\vp}}{2k+1-\delta}\left(|p|^{2k + 1 -\delta} - \frac{1+\delta}{2k+2}\right),
\ee
we thus have shown that
\ber \nonumber 
&&\hspace{-12mm}\frac{1}{2k+1}\io |p|^{2k + 1 -\delta}(x,t) \dd x
+ \int_0^T\io |p|^{2k} |\nabla p|^2 \dd x \dd {t'}
+\frac{1}{2k+2}\int_0^T\ipo\alpha(x)|p|^{2k+2}\dd s(x)\dd {t'}\\ \label{bbe10}
&\le& C\max\{E^{2k+2}, J_k^{1/q'}\}.
\eer


\subsection{Moser iterations} \label{mose}

We first recall a technical lemma proved in \cite[Lemma~3.1]{kg}.

\begin{lemma}\label{l1}
Let $\Omega \subset \real^N$ be a bounded Lipschitzian domain, $N\ge 2$. 
Moreover, let $q_0 = (N+2)/2$, $q_0' = (N+2)/N$, and suppose that the real numbers
 $s, r$  satisfy the inequalities
\be{m1}
\frac12 \le s \le r \le \frac{N+2s}{N+2} \le 1\,.
\ee
Furthermore, assume that a function $v \in L^2(0,T;W^{1,2}(\Omega))$ satisfies for a.~e. 
$t \in (0,T)$ the inequality 
\be{m2}
|v(t)|_{2s}^{2s} + \int_0^T |v({t'})|^2_{W^{1,2}(\Omega)}\dd {t'}
\le A \max\left\{B, \|v\|_{2rq'}\right\}^{2r},
\ee
for some $q' < q_0'$, $A\ge 1$, and $B\ge 1$. Then there exists a constant $C\ge 1$,
which is independent of the choice of $v$, $B$, and $A$, such that
\be{m3}
\|v\|_{2rq_0'} \le C A^{1/(2r)} \max \left\{B,\|v\|_{2rq'}\right\}.
\ee
\end{lemma}

\vspace{2mm}
We now apply Lemma \ref{l1} to the inequality \eqref{bbe10} with $q$ 
given by \eqref{c9}.
Put $v_k := p |p|^{k}$. Then \eqref{bbe10} can be rewritten, using H\"older's inequality, as
\be{m12}
|v_k(t)|_{2s}^{2s} + \int_0^T |v_k({t'})|^2_{W^{1,2}(\Omega)}\dd {t'}
\le (k+1)^2 A \max\left\{E^{k+1}, \|v_k\|_{2rq'}\right\}^{2r}
\ee
with
$$
2s = \frac{2k + 1 - \delta}{k+1}\,, \quad 2r = \frac{2k + 1}{k+1}, \quad q' = \frac{q}{q-1}\,,
$$
and with a constant $A\ge 1$ depending only on the initial and boundary data.
We see that the hypothesis \eqref{m1} of Lemma \ref{l1} is fulfilled whenever $k \ge \delta$.
The assertion of Lemma \ref{l1} then ensures that
\be{m14}
\|v_k\|_{2rq_0'} \le C ((k+1)^2 A)^{1/(2r)} \max \left\{E^{k+1},\|v_k\|_{2rq'}\right\},
\ee
which entails that
\be{m15}
\max\left\{E,\|p\|_{(2k+1)q_0'}\right\} \le C^{1/(k+1)} ((k+1)^2 A)^{1/(2k+1)}\max\left\{E,\|p\|_{(2k+1)q'}\right\}.
\ee
By induction, we check that the choice $b(p) = |p|^{2k}$ is justified for every $k \ge \delta$.
Moreover, we set $\widetilde \nu := (q_0'/q') - 1 > 0$ and define the sequence
$\{k_j\}_{j\ge 0}$ by the formula
\be{m16}
2k_j + 1 = (2\delta +1)(1+\widetilde\nu)^j.
\ee
Set $D_j := \max\left\{E,\|p\|_{(2k_j+1)q_0'}\right\}$. Then \eqref{m15} takes the form
\be{m17}
D_j \le C^{1/(k_j+1)} ((k_j+1)^2 A)^{1/(2k_j+1)} D_{j-1} \ \for \ j\in \nat, 
\ee
and therefore,
\be{m18}
\log D_j - \log D_{j-1} \le \frac{1}{k_j+1} \log C + \frac{1}{2k_j+1} \log((k_j+1)^2 A).
\ee
We have $k_0 = \delta$ and
$D_0 \le C$, by \eqref{c8}--\eqref{c9} and the condition $\delta<1/4$ in 
Hypothesis \ref{h1}\,(iv).
The series on the right-hand side of \eqref{m18} is convergent, and we thus have
$$
D_j \le D_0 \prod_{j=1}^\infty C^{1/(k_j+1)^2}((k_j+1)^2 A)^{1/(2k_j+1)} \le C^*
$$
with a constant $C^*$ independent of $j$, which enables us to conclude that
\be{m19}
\|p\|_\infty \le C^*.
\ee


\subsection{Higher order estimates for the capillary pressure}\label{high}

We aim at taking the limit as  $R\nearrow \infty$ in \eqref{cu1}--\eqref{le7a}. Hence, we can restrict ourselves
to parameter values $R > C^*$ with $C^*$ from \eqref{m19} and rewrite \eqref{cu1}--\eqref{le7a}
in the form
\ber\label{hu1}
\io\left((\dch (\vp(p)+w))_t \eta +\frac{1}{\rho_L}\mu(p)\nabla p\cdot\nabla \eta\right)\dd x
&=& \ipo \alpha(x)(p^*-p)\eta \dd s(x),\qquad \\ \label{hu2}
\nu w_t + \lambda_M w  - p(\chi + \rho^*(1{-}\chi))
- \beta(\hate-\theta_c) &=& - G + H_R(t)\quad a.~e.,\\ \label{hu4}
 \gamma(\hate) \chi_t + \partial I(\chi) - (1-\rho^*)(\Phi(p)+pw)
&\ni& L\left(\frac{\hate}{\theta_c}-1\right)\quad a.~e., \\ \nonumber
\io\left(c_0\theta_t \zeta + \kappa(\hate) \nabla \theta\cdot\nabla \zeta\right)\dd x 
+ \ipo \omega(x)(\theta-\theta^*)\zeta \dd s(x)
&=& \frac{1}{\rho_L} \io \mu(p)Q_R(|\nabla p|^2)\zeta\dd x\\ \label{hu5}
&&\hspace{-110mm}
+\,\io\Big(\chi_t\big((1{-}\rho^*)(\Phi(p) + pw)-L\big)
+ w_t(\dch p {-} \lambda_M w {-} \beta\theta_c {-} G {+} H_R(t))\Big)\zeta\dd x\quad
\eer
for every test functions $\eta,\zeta \in W^{1,2}(\Omega)$, with $\hate = Q_R(\theta)$ and
\be{hu6}
H_R(t) = - \frac{1}{|\Omega|}\io (p \dch + \beta(\hate -\theta_c) -G)(x,t)\dd x\,.
\ee
We test \eqref{hu1} by $\eta = \mu(p) p_t$, which is an admissible choice by \eqref{ge2a}. Then
\ber\nonumber
&&\hspace{-12mm}\io \dch \vp'(p)\mu(p)|p_t|^2 \dd x
+\frac{1}{2\rho_L}\frac{\dd}{\dd t} \io\mu^2(p)|\nabla p|^2\dd x
+ \ipo \alpha(x)(p-p^*)\mu(p) p_t \dd s(x)\\ \label{he1}
&=& \io \left((1-\rho^*)\chi_t w + \dch w_t\right)\mu(p) p_t\dd x.
\eer
Note that by Hypothesis \ref{h1}\,(iv) and \eqref{m19}, we have
$$
\vp'(p) \ge \frac{c_{\vp}}{\max\{1, C^*\}^{1+\delta}}.
$$
We set
\be{mu}
\hat\mu(p) = \int_0^p \tau \mu(\tau) \dd \tau, \quad M(p) = \int_0^p \mu(\tau)\dd\tau,
\ee
and integrate \eqref{he1} in time to obtain the estimate
\ber\nonumber
&&\hspace{-12mm}\int_0^{t}\io |p_t|^2 \dd x\dd {t'} +
\io|\nabla p|^2(x,t)\dd x + \ipo \alpha(x)\hat\mu(p)(x,t) \dd s(x)\\ \nonumber
&\le& C\Bigg({1\,+}\ipo\alpha(x)M(p)|p^*|(x,t) \dd s(x)
+ \int_0^{t} \ipo\alpha(x)M(p)|p^*_t|(x,t') \dd s(x)\dd t'\\ \nonumber
&&+\, \int_0^{t}\io (|\chi_t w| + |w_t|)|p_t|\dd x \dd{t'}\Bigg)\\ \label{he2}
&\le& C\left(1 + \int_0^{t}\io (|\chi_t w| + |w_t|)|p_t|\dd x \dd{t'}\right)
\eer
for all $t \in [0,T]$, whence we infer that
\ber\nonumber
&&\hspace{-12mm}\int_0^{t}\io |p_t|^2 \dd x\dd {t'} +
\io|\nabla p|^2(x,t)\dd x + \ipo \alpha(x)\hat\mu(p)(x,t) \dd s(x)\\ \label{he3}
&\le& C\left(1 + \int_0^{t}\io (|\chi_t w|^2 + |w_t|^2)\dd x \dd{t'}\right).
\eer
By virtue of \eqref{c3}--\eqref{c3a}, \eqref{hu4}, and \eqref{m19}, we have the pointwise bounds
\ber \label{he4}
|w(x,t)| &\le& C \left(1+ \int_0^t \hate(x,{t'})\dd{t'}\right),\\ \label{he5}
|\chi_t(x,t)| &\le& C(1+ |w(x,t)|) \,\le\, C \left(1+ \int_0^t \hate(x,{t'})\dd{t'}\right),\\ \label{he6}
|w_t(x,t)| &\le& C \left(1+ \hate(x,t) + \int_0^t \hate(x,{t'})\dd{t'}\right).
\eer
By \eqref{de2} and the Sobolev embedding theorem, we know that $\hate$ is bounded in
$L^{8/3}(\Omega\times (0,T)) \cap L^2(0,T; L^6(\Omega))$. Let us recall again the Minkowski inequality
$$
\left(\io \left(\int_0^t \hate (x,{t'}) \dd{t'}\right)^6\dd x\right)^{1/6}
\le \int_0^t \left(\io \hate^6(x,{t'})\dd x\right)^{1/6} \dd{t'}\,,
$$
which implies that
\ber \label{he4a}
\io \left(|w(x,t)|^6 + |\chi_t(x,t)|^6\right) \dd x &\le& C \quad \mbox{ for a.\,e. } t \in (0,T),\\ \label{he5a}
\|w_t\|_{8/3} &\le& C.
\eer
Hence, the right-hand side of \eqref{he3} is bounded independently of $R$, and we have for all $t \in [0,T]$ that
\be{he7}
\int_0^t\io |p_t|^2 \dd x\dd {t'} +
\io|\nabla p|^2(x,t)\dd x + \ipo \alpha(x)\hat\mu(p)(x,t) \dd s(x) \le C\,.
\ee
Now let $M(p)$ be as in \eqref{mu}. By \eqref{he4a}--\eqref{he7}, and by comparison in \eqref{hu1},
the term $\Delta M(p)$ is bounded in $L^{2}(\Omega\times (0,T))$, independently of $R$.
In terms of the new variable $\tilde p = M(p)$, the boundary condition \eqref{be2} is nonlinear,
and the $W^{2,2}$-regularity of $ M(p)$ follows from considerations similar to those used
in the proof of \cite[Theorem 4.1]{kpa}, 
inspired by \cite{jn}.  We thus may employ the Gagliardo-Nirenberg inequality \eqref{gnm} in the form
\be{gnm}
|\nabla M(p)(t)|_q \le C\left(|\nabla M(p)(t)|_2+|\nabla M(p)(t)|_2^{1-\rho}|\Delta M(p)(t)|_2^\rho\right)
\ee
with $\rho = 3(\frac12 - \frac{1}{q})$. 
Together with \eqref{he7}, we conclude that
\be{he8a}
\int_0^T|\nabla p(t)|_q^s \dd t \le C \quad \for q \in (2,6] \ \mbox{ and }\ \frac{1}{q} + \frac{2}{3s} = \frac12.
\ee
In particular, for $s=4$ and $s=q$ we obtain, respectively,
\be{he8}
\int_0^T|\nabla p(t)|_3^4 \dd t \le C\,, \quad \|\nabla p\|_{10/3} \le C.
\ee


\subsection{Higher order estimates for the temperature}\label{temp}

The previous estimates \eqref{he4a}--\eqref{he5a}
and \eqref{he8} entail that \eqref{hu5} has the form
\be{te1}
\io\left(c_0\theta_t \zeta
+ \kappa(\hate) \nabla \theta\cdot\nabla \zeta\right)\dd x 
+ \ipo \omega(x)(\theta-\theta^*)\zeta \dd s(x) = \io \tilde F \zeta\dd x
\ee
for every $\zeta \in W^{1,2}(\Omega)$, with a function $\tilde F$ such that
\be{te2}
\|\tilde F\|_{5/3} \le C\,, \quad \int_0^T|\tilde F(t)|_{3/2}^2 \dd t \le C\,,
\ee
independently of $R$.

Assume now that for some $p_0 \ge 8/3$ we have
\be{te3}
\|\hate\|_{p_0} \le C.
\ee
We know that this is true for $p_0 = 8/3$ by virtue of \eqref{de2}. Set $r_0 = 2p_0/5$.
Then we may put $\zeta = \hate^{r_0}$ in \eqref{te1} and obtain, using Hypothesis \ref{h1}\,(ii), that
\be{te4}
\frac{1}{r_0+1}\io \hate^{r_0 + 1}(x,t)\dd x + r_0 \int_0^t\io \hate^{r_0 + a} |\nabla\hate|^2\dd x\dd{t'}
\le C\,.
\ee
We now denote
$$
v = \hate^p\,, \quad p = 1 + \frac{r_0+a}{2}\,, \quad s = \frac{r_0 + 1}{p}\,,
$$
and rewrite \eqref{te4} as
\be{te5}
\io |v|^{s}(x,t)\dd x + \int_0^t\io |\nabla v|^2\dd x\dd{t'} \le C(r_0 + 1)\,.
\ee
By the Gagliardo-Nirenberg inequality \eqref{gn}, we have $\|v\|_q \le C(r_0 + 1)$ for $q = 2 + \frac{2s}{3}$.
Hence,
\be{te6}
\|\hate\|_{p_1} \le C(r_0 + 1) \quad \for \ p_1 = pq = \frac{2p_0}{3} + \frac83 + a\,.
\ee
We now proceed by induction according to the recipe
$
p_{j+1} = \frac{2p_j}{3} + \frac83 + a\,,\  r_j = \frac{2p_j}{5}\,.
$
We have $\lim_{j\to \infty} p_j = 8 + 3a$. After finitely many steps, we
may stop the algorithm and put $\bar p := p_j < 8 + 3a$ with
\be{te7}
\|\hate\|_{\bar p} + \supess |\hate(t)|_{\bar r+ 1} \le C\,,\quad
\bar r = \frac{2\bar p}{5} > \hat a,
\ee
with the constant  $\hat a$ introduced in Hypothesis \ref{h1}\,(ii).
By Proposition \ref{t2}, we may test \eqref{te1} by $\theta$, which yields
\be{te8}
\io \theta^{2}(x,t)\dd x + \int_0^t\io \kappa(\hate) |\nabla\theta|^2\dd x\dd{t'}
+ \int_0^t\ipo \omega(x) \theta^2 \dd s(x)\dd{t'} \le C \|\theta\|_{5/2}\,.
\ee
Using the Gagliardo-Nirenberg inequality again, for instance, we conclude that
\be{te9}
\io \theta^{2}(x,t)\dd x + \int_0^t\io \kappa(\hate) |\nabla\theta|^2\dd x\dd{t'}
+ \int_0^t\ipo \omega(x) \theta^2 \dd s(x)\dd{t'} \le C\,.
\ee
This enables us to derive an upper bound for the integral
$\io \kappa(\hate) \nabla\theta\cdot\nabla\zeta \dd x$,
which we need for getting an estimate for $\theta_t$ from the equation \eqref{te9}. We have,
by H\"older's inequality and Hypothesis \ref{h1}\,(ii), that
\ber \nonumber
\io |\kappa(\hate) \nabla\theta\cdot\nabla\zeta| \dd x &=&
\io |\kappa^{1/2}(\hate) \nabla\theta\cdot\kappa^{1/2}(\hate)\nabla\zeta| \dd x\\
\label{e715}
&\le& C\left(\io\kappa(\hate) |\nabla\theta|^2\dd x\right)^{1/2}
 \left(\io \hate^{1+\hat a}|\nabla\zeta|^2\dd x\right)^{1/2}.
\eer
We now choose $\hat q > 1$ such that $(1+\hat a)\hat q = 1+\bar r$, where $\bar r$ 
is defined in \eqref{te7}.
Choosing now
\be{qstar}
q^* = \frac{2 \hat q}{\hat q - 1},
\ee
we obtain from H\"older's inequality that
\be{e714}
\io \hat\theta^{1+\hat a}|\nabla\zeta|^2\dd x \le \left(\io \hat\theta^{1+\bar r} \dd x\right)^{1/\hat q}
\left(\io |\nabla\zeta|^{q^*}\dd x\right)^{2/q^*}
\le C \left(\io |\nabla\zeta|^{q^*}\dd x\right)^{2/q^*}\,,
\ee
by virtue of \eqref{te7}. Eq.~\eqref{e715} then yields the bound
\be{e716}
\io |\kappa(\hat\theta) \nabla\theta\cdot\nabla\zeta| \dd x 
\le C\left(\io\kappa(\hate) |\nabla\theta|^2\dd x\right)^{1/2}
\left(\io |\nabla\zeta|^{q^*}\dd x\right)^{1/q^*}.
\ee
Hence, by \eqref{te9},
\be{e717}
\int_0^T\io |\kappa(\hate) \nabla\theta\cdot\nabla\zeta| \dd x \dd t
\le C \|\zeta\|_{L^2(0,T;W^{1,q^*}(\Omega))}\,.
\ee
{}From \eqref{te2} it follows that testing with $\zeta \in L^2(0,T;W^{1,q^*}(\Omega))$ is admissible.
We thus obtain from \eqref{te1} that
\be{e718}
\int_0^T\io \theta_t \zeta \dd x \dd t \le C \|\zeta\|_{L^2(0,T;W^{1,q^*}(\Omega))}\,.
\ee


\section{Proof of Theorem \ref{t1}} \label{proo}

Let $R_i \nearrow \infty$ be a sequence such that $R_1 > C^*$ with $C^*$ as in \eqref{m19}, and let
$(p,w,\chi,\theta) = (p\oi, w\oi, \chi\oi, \theta\oi)$ be
solutions of \eqref{hu1}--\eqref{hu6} corresponding to $R = R_i$,
with $\hate = \hate\oi = Q_{R_i}(\theta\oi)$ and test functions $\eta,\zeta \in W^{1,2}(\Omega)$.
Our aim is to check that at least a subsequence converges as $i \to \infty$
to a solution of \eqref{wu1}--\eqref{wu5}, \eqref{le7},
with test functions $\eta\in W^{1,2}(\Omega)$, $\zeta \in W^{1,q^*}(\Omega)$ with $q^*$ as in Theorem \ref{t1}.

First, for the capillary pressure $p = p\oi$ we have the estimates \eqref{m19}, \eqref{he7}, \eqref{he8},
which imply that, passing to a subsequence if necessary,
\bers
p\oi \to p && \mbox{strongly in } L^r(\Omega\times (0,T)) \ \mbox{ for every } \ r\ge 1\,,\\
p_t\oi \to p_t && \mbox{weakly in } L^2(\Omega\times (0,T))\,,\\
\nabla p\oi \to \nabla p && \mbox{strongly in } L^r(\Omega\times (0,T)) \ \mbox{ for every }
\ 1 \le r < \frac{10}{3}\,.
\eers
We easily show that
\be{pr2}
Q_{R_i}\bigl(|\nabla p\oi|^2\bigr) \to |\nabla p|^2 \ \mbox{ strongly in } L^r(\Omega\times (0,T))
\ \mbox{ for every } \ 1 \le r < \frac{5}{3}\,.
\ee
Indeed, let $\Omega\oi_T \subset \Omega\times (0,T)$ be the set of all $(x,t) \in \Omega\times (0,T)$
such that $|\nabla p\oi(x,t)|^2 > R_i$. By \eqref{he8}, we have
$$
C \ge \int_0^T\io |\nabla p\oi(x,t)|^{10/3} \dd x\dd t \ge \iint_{\Omega\oi_T}|\nabla p\oi(x,t)|^{10/3} \dd x\dd t
\ge |\Omega\oi_T| R_i^{5/3}\,,
$$ 
hence $|\Omega\oi_T| \le C R_i^{-5/3}$. For $r < \frac{5}{3}$, we use H\"older's inequality to get the estimate
\bers
&&\hspace{-12mm}\int_0^T\io \left|Q_{R_i}(|\nabla p\oi|^2) - |\nabla p\oi|^2\right|^r\dd x\dd t
= \iint_{\Omega\oi_T}\left|R_i - |\nabla p\oi|^2\right|^r\dd x\dd t
\le \iint_{\Omega\oi_T} |\nabla p\oi|^{2r}\dd x\dd t \\
&\le& \left(\iint_{\Omega\oi_T} |\nabla p\oi|^{10/3}\dd x\dd t\right)^{3r/5} |\Omega\oi_T|^{1-(3r/5)}\,,
\eers
and \eqref{pr2} follows.

For the temperature $\theta = \theta\oi$, we proceed in a similar way. From the compactness result
in \cite[Theorem 5.1]{li}, it follows that, for a subsequence,
$$
\theta\oi \to \theta \ \mbox{ strongly in } L^2(\Omega\times (0,T)).
$$
Furthermore, by \eqref{te7}, $\hate\oi$ are uniformly bounded in $L^{r}(\Omega\times (0,T))$
for every $r < 8+3a$. A similar argument as above yields that
$$
\hate\oi \to \theta \ \mbox{ strongly in } L^r(\Omega\times (0,T)) \ \mbox{ for every }
\ 1 \le r < 8+3a\,.
$$
Indeed, by \eqref{te9} and \eqref{e718},
\bers
\theta_t\oi \to \theta_t && \mbox{weakly in } L^2(0,T; W^{-1,q^*}(\Omega))\,,\\
\nabla \theta\oi \to \nabla \theta && \mbox{weakly in } L^2(\Omega\times (0,T))\,.
\eers
The strong convergences of $w\oi \to w$, $w\oi_t \to w_t$, $\chi\oi \to \chi$, $\chi\oi_t \to \chi_t$
are handled using the estimates \eqref{he4}--\eqref{he6}
similarly as in the proof of Proposition \ref{t2} at the end of Section \ref{cut}.
This enables us to pass to the limit as $R\nearrow \infty$ in the system \eqref{hu1}--\eqref{hu6}
and thus to complete the proof of Theorem \ref{t1}.


{\small
\begin{thebibliography}{99}

\bibitem{albers} B. Albers:\ \emph{Modeling and Numerical Analysis of Wave
Propagation in Saturated and Partially Saturated Porous Media}. Habilitation
Thesis. Ver\"{o}ffentlichungen des Grundbauinstitutes der Technischen
Universit\"{a}t Berlin, Shaker: Aachen; 2010.


\bibitem{ak} B. Albers and P. Krej\v{c}\'{\i}:\
Unsaturated porous media flow with thermomechanical interaction.
{\em Math. Meth. Appl. Sci.\/} {\bf 39} (2016),  2220--2238.

\bibitem{bin} O.~V. Besov, V.~P. Il'in, and S.~M. Nikol'ski\u{\i}:\
\emph{Integral Representations of Functions and Imbedding Theorems}.
Scripta Series in Mathematics. Halsted Press (John Wiley \& Sons):
New York-Toronto, Ont.-London; 1978 (Vol. I), 1979 (Vol. II).
Russian version Nauka: Moscow; 1975.

\bibitem{bs} 
  M. Brokate and J. Sprekels:\
  \emph{Hysteresis and Phase Transitions.}
  Applied Mathematical Sciences Vol. 121,
  Springer-Verlag,
  New York (1996).

\bibitem{dkr1} B. Detmann, P. Krej\v c\'\i, and E. Rocca:\
Solvability of an unsaturated porous media flow problem with thermomechanical interaction.
{\em SIAM J. Math. Anal.\/} {\bf 48} (2016), 4175--4201.


\bibitem{dkr2} B.~Detmann, P.~Krej\v c\'\i, and E.~Rocca:\
Periodic waves in unsaturated porous media with hysteresis, 
preprint arXiv:1606.04665, 1--15, to appear in the Proceedings of ECM 2016.

\bibitem{fremond} 
M. Fr\'emond:\ \emph{Non-Smooth Thermomechanics.} Springer-Verlag, Berlin (2002).

\bibitem{fr1} 
M.~Fr\'emond and E.~Rocca:\ Well-posedness of a phase transition model
with the possibility of voids, {\em Math. Models Methods Appl. Sci.}
{\bf 16} (2006), 559--586.

\bibitem{fr2} 
M.~Fr\'emond and E.~Rocca:\ Solid liquid phase changes with
different densities, \emph{Q. Appl. Math.} {\bf 66} (2008), 609--632.

\bibitem{cmuc} P. Krej\v c\'\i:\
Hysteresis operators -- a new approach to evolution differential inequalities.
{\em Comment. Math. Univ. Carolinae} {\bf 33} (1989), 525--536.

 \bibitem{mb} P. Krej\v c\'\i:\
Elastoplastic reaction of a container to water freezing.
{\em Math. Bohem.} {\bf 135} (2010), 423--441.

\bibitem{kg} P. Krej\v c\'\i:\ Boundedness of solutions to a degenerate diffusion equation.
To appear in: {\em Solvability, Regularity, Optimal Control of Boundary Value Problems 
for PDEs} (P. Colli, A. Favini, E. Rocca, G. Schimperna, J. Sprekels, eds.), Springer INdAM Series. 

\bibitem{kpa} P.~Krej\v{c}\'{\i} and L. Panizzi:\
Regularity and uniqueness in quasilinear parabolic systems.
{\em Appl. Math.} {\bf 56} (2011), 341--370.

\bibitem{kr} P. Krej\v{c}\'{\i} and E. Rocca:\ 
Well-posedness of an extended model for water-ice phase transitions.
\emph{Discrete Contin. Dyn. Syst. Ser. S} {\bf 6} (2013), 439--460. 

\bibitem{krsbottle} 
P.~Krej\v c\'{\i}, E.~Rocca, and J.~Sprekels:\ A bottle in a
freezer. \emph{SIAM J. Math. Anal.}  {\bf 41}  (2009), 1851--1873.

\bibitem{krsgrav} 
P. Krej\v c\'{\i}, E. Rocca, and J.  Sprekels:\
Phase separation in a gravity field. \emph{Discrete Contin. Dyn. Syst. Ser.
S} {\bf 4} (2011), 391--407.

\bibitem{krsWil}
P. Krej\v c\'{\i}, E. Rocca, and J.  Sprekels:\
Liquid-solid phase transitions in a deformable container. 
Contribution to the book {\em Continuous Media with Microstructure}
on the occasion of Krzysztof Wilmanski's 70th birthday, Springer
(2010), 285--300.


\bibitem{kss2} 
{P.~Krej\v{c}\'{\i}, J.~Sprekels, and U.~Stefanelli}:\
{Phase-field models with hysteresis in one-dimensional
thermoviscoplasticity}. \emph{SIAM J. Math. Anal.} {\bf 34} (2002), 409--434.

\bibitem{kss} 
{P.~Krej\v{c}\'{\i}, J.~Sprekels, and U.~Stefanelli}:\
{One-dimensional thermo-visco-plastic processes with hysteresis and
phase transitions}. {\em Adv. Math. Sci. Appl.} {\bf 13} (2003), 695--712.

\bibitem{li} J.-L. Lions:\ \emph{Quelques m\'ethodes de r\'esolution des
probl\`emes aux limites non lin\'eaires.} Dunod; Gauthier-Villars, Paris 1969.

\bibitem{jn}
J.~Ne\v cas:\ \emph{Les m\'ethodes directes en th\'eorie des \'equations elliptiques}.
Academia, Prague, 1967.

\bibitem{rr1} 
{E.~Rocca and R.~Rossi}:\
{A nonlinear degenerating PDE system modelling phase transitions
in thermoviscoelastic materials}. {\em J. Differential Equations} {\bf 245} (2008), 3327--3375.

\bibitem{rr2} 
{E.~Rocca and R.~Rossi}:\
{ Global existence of strong solutions to the one-dimensional full
model for phase transitions in thermoviscoelastic materials}.
{\em Appl. Math.} {\bf 53} (2008), 485--520.

\bibitem{RocRos12}
	E. Rocca and R. Rossi:\ A degenerating PDE system for
	phase transitions and damage. {\em  Math. Models Methods Appl. Sci.} {\bf 24} (2014), 1265--1341.

\bibitem{RocRos14}
E. Rocca and R. Rossi:\
``Entropic'' solutions to a thermodynamically consistent PDE system for phase transitions and damage.
{\em SIAM J. Math. Anal.} {\bf  47} (2015), 2519--2586. 

\bibitem{show} R.~E. Showalter:\
Diffusion in deforming porous media.
\emph{Dyn. Contin. Discrete Impuls. Syst. Ser. A Math. Anal.} {\bf 10} (2003), 661--678.

\bibitem{ss} R.~E. Showalter and U. Stefanelli:\
Diffusion in poro-plastic media.
\emph{Math. Methods Appl. Sci.} {\bf 27} (2004), 2131--2151. 

\bibitem{wil}
K. Wilma{\'n}ski:\
Macroscopic modeling of porous materials. In: Mathematical Modelling and Analysis in Continuum Mechanics of Microstructured Media:
Professor Margaret Wo{\' z}niak pro memoria: sapiens mortem non timet. Awrejcewicz J, et al. (eds). Wydawnictwo Politechniki \'Slaskiej: Gliwice, 2010;
167–195.

\bibitem{visintin} 
A.~Visintin:\
\emph{Models of Phase Transitions}.,
Progress in Nonlinear Differential Equations and their Applications Vol. 28,
Birkh\"auser, Boston, 1996.


\end{thebibliography}
}
\end{document}